\def\yes{\if00}
\def\no{\if01}
\def\iftwelvept{\yes}
\def\ifusepdf{\no}
\def\ifusepsfont{\no}

\iftwelvept
\documentclass[leqno,12pt]{amsart}
\else
\documentclass[leqno]{amsart}
\fi
\usepackage{amsfonts}
\usepackage{amssymb}
\usepackage{amscd}
\usepackage{latexsym}
\usepackage{verbatim}
\usepackage[all]{xy}
\ifusepdf
\usepackage{hyperref}
\else\fi
\ifusepsfont
\usepackage[T1]{fontenc}
\usepackage{times}
\else\fi

\iftwelvept
\else\fi

\newtheorem{Theorem}{Theorem}[section]
\newtheorem{Proposition}[Theorem]{Proposition}
\newtheorem{Lemma}[Theorem]{Lemma}
\newtheorem{Corollary}[Theorem]{Corollary}
\newtheorem{Claim}{Claim}[Theorem]

\theoremstyle{definition}

\newtheorem{Remark}[Theorem]{Remark}

\renewcommand{\theTheorem}{\arabic{section}.\arabic{Theorem}}

\newcommand{\ZZ}{{\mathbb{Z}}}
\newcommand{\QQ}{{\mathbb{Q}}}
\newcommand{\RR}{{\mathbb{R}}}
\newcommand{\CC}{{\mathbb{C}}}
\newcommand{\PP}{{\mathbb{P}}}
\newcommand{\OO}{{\mathcal{O}}}
\newcommand{\Aff}{{\mathbb{A}}}

\newcommand{\Spec}{\operatorname{Spec}}
\newcommand{\Pic}{\operatorname{Pic}}
\newcommand{\ord}{\operatorname{ord}}

\newcommand{\Supp}{\operatorname{Supp}}

\newcommand{\Coker}{\operatorname{Coker}}

\newcommand{\Proof}{{\sl Proof.}\quad}

\newcommand{\QED}{{\unskip\nobreak\hfil\penalty50\quad\null\nobreak\hfil
{$\Box$}\parfillskip0pt\finalhyphendemerits0\par\medskip}}
\newcommand{\rest}[2]{\left.{#1}\right\vert_{{#2}}}
\begin{document}

%\begin{flushright}
%{\fbox{Second Draft of math.NT/0405007}
%\footnote{
%I would like to replace math.NT/0405007 by 
%this draft as the second version. Compared with the previous draft, 
%we treat polynomial automorphisms of dynamical degree $\geq 2$, 
%not just regular polynomial automorphisms. 
%Also the proof of (0.2) is simplified thanks 
%to Prof. Noboru Nakayama.}}
%\end{flushright}
%\bigskip

%%%%%%%%%%%
%% Title %%
%%%%%%%%%%%
\title[Canonical height functions]%
{Canonical height functions on 
the affine plane associated with polynomial automorphisms}
\author{Shu Kawaguchi}
\address{Department of Mathematics, Faculty of Science,
Kyoto University, Kyoto, 606-8502, Japan}
\email{kawaguch@math.kyoto-u.ac.jp}
%\date{\DateTime, (\Version)}
%\thanks{{\em 2000 Mathematics Subject Classification:} 11G50}
\begin{abstract}
Let $f: \Aff^2\to\Aff^2$ be a polynomial automorphism 
of dynamical degree $\delta \geq 2$ over a number field $K$. 
(This is equivalent to say that $f$ is a polynomial automorphism 
that is not triangularizable.) Then we construct 
canonical height functions defined on $\Aff^2(\overline{K})$ 
associated with $f$. 
These functions satisfy the Northcott finiteness property, 
and an $\overline{K}$-valued point on $\Aff^2(\overline{K})$ 
is $f$-periodic if and only if its height is zero. 
As an application of canonical height functions, 
we give an estimate on the number of points with bounded height 
in an infinite $f$-orbit. 
\end{abstract} 

\maketitle

\section*{Introduction and the statement of the main results}
\renewcommand{\theTheorem}{\Alph{Theorem}} 
One of the basic tools in Diophantine geometry is the theory of 
height functions. On Abelian varieties defined over a number field, 
N\'eron and Tate developed the theory of canonical height functions 
that behave well relative to the $[n]$-th power 
map (cf. \cite[Chap.~5]{La}). 
On certain K3 surfaces with two involutions, 
Silverman \cite{SiK3} developed the theory of canonical height 
functions that behave well relative to the two involutions. 
For the theory of canonical height functions on some other 
projective varieties, see for example \cite{CS}, \cite{Zh}, 
\cite{Ka}. In this paper, we show the existence of canonical height
functions on the affine plane associated with polynomial automorphisms of dynamical degree $\geq 2$.

Consider a polynomial automorphism $f: \Aff^2 \to \Aff^2$ 
given by
\[
f \begin{pmatrix} x \\ y \end{pmatrix} 
= \begin{pmatrix} p(x,y)  \\ q(x,y) \end{pmatrix},  
\]
where $p(x,y)$ and $q(x,y)$ are polynomials in two variables. 
The degree $d$ of $f$ is defined by $d := \max\{\deg p, \deg q\}$. 
The dynamical degree $\delta$ of $f$ is defined by 
\[
\delta := \lim_{n\to +\infty} \left(\deg f^n\right)^{\frac{1}{n}},  
\] 
which is an integer with $1 \leq \delta \leq d$. 
We let $d \geq 2$. 

Polynomial automorphisms with $\delta = d$ are exactly regular
polynomial automorphisms. Here a polynomial automorphism $f: \Aff^2
\to \Aff^2$ is said to be regular if the unique point of indeterminacy
of $\overline{f}$ is different from the unique point of 
indeterminacy of $\overline{f^{-1}}$, where the birational map
$\overline{f}: \PP^2 \dasharrow \PP^2$ (resp. $\overline{f^{-1}}: \PP^2
\dasharrow \PP^2$) is the extension of $f$ (resp. $f^{-1}$).  
In the moduli of polynomial automorphisms of degree $d$, 
regular polynomial automorphisms constitute general members, 
including H\'enon maps.  

The other extreme is polynomial automorphisms of dynamical degree
$\delta = 1$,
and they are exactly triangularizable automorphisms. Here a polynomial
automorphism $f: \Aff^2 \to \Aff^2$ is said to be triangularizable if it
is conjugate, in the group of polynomial automorphisms, to a
polynomial automorphism of the form
\[
f \begin{pmatrix} x \\ y \end{pmatrix} 
= \begin{pmatrix} a x + P(y) \\ b y + c  \end{pmatrix}, 
\]
where $ab \neq 0$ and $P(y)$ is a polynomial in $y$.  
For more details, see the survey 
of Sibony \cite{Sib} and the references therein. 
See also \S\ref{sec:geometry}.

Over a number field, Silverman \cite{SiHenon} studied 
arithmetic properties of quadratic H\'enon maps,  
and then Denis \cite{De} 
studied arithmetic properties of H\'enon maps 
and some classes of polynomial automorphisms. 
Marcello \cite{Ma1}, \cite{Ma}
studied arithmetic properties of some other classes of polynomial 
automorphisms of the affine spaces, 
including regular polynomial automorphisms. 

Our first result shows the existence of 
height functions that behave well relative to 
polynomial automorphisms of $\Aff^2$. 

\begin{Theorem}
\label{thm:main}
Let $f: \Aff^2 \to \Aff^2$ be a polynomial automorphism of 
dynamical degree $\delta \geq 2$ over a number field $K$. 
\textup{(}This is equivalent to say 
that $f$ is a polynomial automorphism 
that is not triangularizable.\textup{)} 
Then there exists a function 
$\widehat{h}: \Aff^2(\overline{K}) \to \RR$
with the following properties: 
\begin{enumerate}
\item[(i)] 
$h_{nv} \gg\ll \widehat{h}$ on $\Aff^2(\overline{K})$
\textup{(}Here $h_{nv}$ is the logarithmic naive height function, and 
$h_{nv} \gg\ll \widehat{h}$ means that there are positive constants 
$a_1, a_2$ and constants $b_1, b_2$ such that 
$a_1 h_{nv} + b_1 \leq \widehat{h} \leq a_2 h_{nv} + b_2$\textup{)}
\textup{;} 
\item[(ii)] 
$\widehat{h} \circ f + \widehat{h} \circ f^{-1}
= \left(\delta + \frac{1}{\delta} \right) \widehat{h}$. 
\end{enumerate}
Moreover, $\widehat{h}$ enjoys the following uniqueness property\textup{:}  
if $\widehat{h}'$ is another function satisfying 
(i) and (ii) 
such that 
$\widehat{h}' = \widehat{h} + O(1)$, then 
$\widehat{h}' = \widehat{h}$. 
We call a function $\widehat{h}$ satisfying 
(i) and (ii)  
a {\em canonical height function} associated with the 
polynomial automorphism $f$. 
\end{Theorem}

It follows from (i) that $\widehat{h}$ satisfies 
the Northcott finiteness property. Namely, 
for any positive number $M$ and positive integer $D$, the set
$\{ x \in \Aff^2(\overline{K}) 
\mid [K(x):K] \leq D, \; \widehat{h}(x) \leq M \}$
is finite. This leads to the following corollary, which shows that 
the set of $\overline{K}$-valued $f$-periodic points is not only 
a set of bounded height but also characterized as the set of 
height zero with respect to a canonical height function
associated with $f$. 

\begin{Corollary}
\label{cor:main}
With the notation and assumption in Theorem~\ref{thm:main},  
\begin{enumerate}
\item[(1)] 
$\widehat{h}(x) \geq 0$ for all $x \in \Aff^2(\overline{K})$. 
\item[(2)] 
$\widehat{h}(x) = 0$ if and only if $x$ is $f$-periodic. 
\textup{(}Here, $x \in \Aff^2(\overline{K})$ is said to be $f$-periodic 
if $f^m(x) = x$ for some positive integer $m$.\textup{)}
\end{enumerate}
\end{Corollary}

As an application of canonical height functions, 
we obtain an estimate on the number of points with bounded 
height in an infinite $f$-orbit. 
First we introduce some notation and terminology.  
For a canonical height function $\widehat{h}$ associated 
with $f$, we set
\[
\widehat{h}^+(x) = \frac{\delta^2}{\delta^4 -1}
\left(\delta \widehat{h}(f(x)) 
- \frac{1}{\delta} \widehat{h}(f^{-1}(x))\right), 
\quad 
\widehat{h}^-(x) = \frac{\delta^2}{\delta^4 -1}
\left(\delta \widehat{h}(f^{-1}(x)) 
- \frac{1}{\delta} \widehat{h}(f(x))\right). 
\]
Then $\widehat{h}^+ \geq 0$ and $\widehat{h}^- \geq 0$, 
and $\widehat{h}^+(x) = 0$ if and only if 
$\widehat{h}^-(x) = 0$ if and only if $x$ is $f$-periodic 
(cf. Lemma~\ref{lemma:h:+:-}).
For a point $x \in \Aff^2(\overline{K})$, 
let $O_f(x) := \{f^l(x) \mid l \in \ZZ\}$ denote the $f$-orbit 
of $x$. For a non $f$-periodic point $x \in \Aff^2(\overline{K})$, 
we set
\[
\widehat{h}(O_f(x)) 
= \frac{\log\left(\widehat{h}^+(y) \widehat{h}^-(y)\right)}{\log\delta}
\]
for any $y \in O_f(x)$. Then $\widehat{h}(O_f(x))$ is 
well-defined, i.e., $\widehat{h}(O_f(x))$ is independent of 
the choice of $y \in O_f(x)$. Moreover, as a function of 
$x$, we have 
$\widehat{h}(O_f(x)) \gg\ll \min_{y \in O_f(x)} \log \widehat{h}(y)$  
on $\Aff^2(\overline{K})\setminus\{\text{$f$-periodic points}\}$ 
(cf. Lemma~\ref{lemma:height:of:orbit}). 

For regular polynomial automorphisms of degree $d \geq 2$, 
it is known that, for a non $f$-periodic point $x \in \Aff^2(\overline{K})$, 
one has $\lim_{T \to +\infty} 
\frac{\#\{ y \in O_f(x) \mid h_{nv}(y) \leq T\}}{\log T} 
= \frac{2}{\log d}$  
(\cite[Theorem~C]{SiHenon}, \cite[Th\'eor\`eme~2]{De}, 
and \cite[Th\'eor\`eme~A]{Ma}).  
The next theorem gives its refinement and generalization. 

\begin{Theorem}
\label{thm:estimate}
Let $f: \Aff^2 \to \Aff^2$ be a polynomial automorphism of 
dynamical degree $\delta \geq 2$ over a number field $K$. 
Suppose $x \in \Aff^2(\overline{K})$ is not an $f$-periodic point. 
Then, 
\begin{equation}
\label{eqn:estimate}
\#\{ y \in O_f(x) \mid h_{nv}(y) \leq T\}
= \frac{2}{\log\delta} \log T - \widehat{h}(O_f(x)) + O(1)
\quad \text{as $T\to +\infty$},
\end{equation} 
where the $O(1)$ constant depends only on $f$ and the choice 
of $\widehat{h}$. 
\end{Theorem}

It seems interesting that the dynamical degree of $f$ appears 
in the left-hand side of \eqref{eqn:estimate}. 
We remark that, when $f$ is not regular, 
i.e., $(2 \leq)\ \delta < \deg f$ , 
even a weaker estimate $\lim_{T \to  +\infty} 
\frac{\#\{ y \in O_f(x) \mid h_{nv}(y) \leq T\}}{\log T} 
= \frac{2}{\log\delta}$ seems new. 

The contents of this paper is as follows. 
In \S\ref{sec:preliminaries} 
we briefly review the properties of height functions. 
In \S\ref{sec:geometry:general} 
we show that if $f$ is a regular polynomial automorphism of degree 
$d \geq 2$ then there is a constant $c$ such that
\begin{equation}
\label{eqn:Henon:blowups}
h_{nv}(f(x)) + h_{nv}(f^{-1}(x)) 
\geq \left(d + \frac{1}{d}\right) h_{nv}(x) - c  
\end{equation}
for all $x \in \Aff^2(\overline{K})$. 
In \S\ref{sec:geometry} we recall H\'enon maps, 
Friedland--Milnor's theorem on the conjugacy classes of 
polynomial automorphisms, and some properties of 
dynamical degrees of polynomial automorphisms. 
In \S\ref{sec:canonical:height:function} we prove Theorem~\ref{thm:main} and 
Corollary~\ref{cor:main} in a more general setting 
of polynomial automorphisms of $\Aff^n$ 
whose conjugates satisfy an inequality similar 
to \eqref{eqn:Henon:blowups}. In \S\ref{sec:estimate} 
we prove Theorem~\ref{thm:estimate} 
in this more general setting. On certain K3 surfaces, Silverman 
counted the number of points with bounded height in a given 
infinite chain (\cite[\S3]{SiK3}). 
Our method of proof of Theorem~\ref{thm:estimate} is inspired 
by his method. 

\smallskip
{\sl Acknowledgments.}\quad 
The author expresses his sincere gratitude to Prof. Noboru Nakayama 
for simplifying the proof of \eqref{eqn:Henon:blowups}. 

\medskip
\renewcommand{\theTheorem}{\arabic{section}.\arabic{Theorem}}
\setcounter{equation}{0}
\section{Quick review on height theory}
\label{sec:preliminaries}
In this section, we briefly review the properties of height 
functions that we will use in this paper. 

Let $K$ be a number field and $O_K$ its ring of integers. 
For $x = (x_0:\cdots:x_n) \in \PP^n(K)$, the logarithmic naive 
height of $x$ is defined by 
\[
h_{nv}(x) 
= \frac{1}{[K:\QQ]}\left[ \sum_{P \in \Spec(O_K)\setminus\{0\}} 
\max_{0 \leq i \leq n} \{ -\ord_P(x_i) \} \log \#(O_K/P)
+ \sum_{\sigma: K \hookrightarrow \CC} 
\max_{0 \leq i \leq n} \{\log \vert\sigma(x_i)\vert \} \right]. 
\] 
This definition naturally extends to all 
points $x \in \PP^n(\overline{\QQ})$ as to give 
the logarithmic naive height function 
$h_{nv}: \PP^n(\overline{\QQ}) \to \RR$. 

We begin by the following two basic properties 
of height functions. 

\begin{Theorem}[Northcott's finiteness theorem, \cite{SiHt} Corollary~3.4]
\label{thm:Northcott}
For any positive number $M$ and positive integer $D$, the set
\[
\left\{ x \in \PP^n(\overline{\QQ}) 
\mid [\QQ(x):\QQ] \leq D, \ h_{nv}(x) \leq M \right\}
\]
is finite. 
\end{Theorem}

\begin{Theorem}[\cite{SiHt} Theorem~3.3, \cite{La} Chap.~4, Prop.~5.2]
\label{thm:height:machine}
\begin{enumerate}
\item[(1)]
\textup{(Height machine)}\; 
There is a unique way to attach, 
for any projective variety $X$ defined over $\overline{\QQ}$, 
a map 
\[
h_X : \Pic(X) \longrightarrow 
\frac{\text{\{real-valued functions on $X(\overline{\QQ})$\}}}%
{\text{\{real-valued bounded functions on $X(\overline{\QQ})$\}}}, 
\quad L \mapsto h_{X,L}
\]
with the following properties\textup{:} 
\begin{enumerate}
\item[(i)] $h_{X, L\otimes M} = h_{X, L} + h_{X, M} + O(1)$ 
for any $L, M \in \Pic(X)$\textup{;} 
\item[(ii)] If $X = \PP^n$ and $L= \OO_{\PP^n}(1)$, then 
$h_{\PP^n, \OO_{\PP^n}(1)} = h_{nv} + O(1)$\textup{;} 
\item[(iii)] If $f: X \to Y$ is a morphism of projective varieties 
and $L$ is a line bundle on $X$, then 
$h_{X, f^*L} = h_{Y, L}\circ f + O(1)$. 
\end{enumerate}
\item[(2)]
\textup{(Positivity of height)}\;
Let $X$ be projective variety defined over $\overline{\QQ}$ 
and $L$ a line bundle on $X$. We set 
$B = \Supp(\Coker(H^0(X,L)\otimes\OO_X \to L))$. 
Then there exists a constant $c_1$ such that 
$h_{X, L}(x) \geq c_1$
for all $x \in (X\setminus B)(\overline{\QQ})$. 
\end{enumerate}
\end{Theorem} 

A rational map $f=[F_0:F_1:\cdots:F_n]: \PP^n \dasharrow \PP^n$ 
defined over $\overline{\QQ}$ is said to be of degree $d$ if 
the $F_i$'s are homogeneous polynomials of degree $d$ 
over $\overline{\QQ}$, with no common factors. 
Let $I_f \subset \PP^n(\overline{\QQ})$ 
denote the locus of indeterminacy. 

\begin{Theorem}[\cite{La} Chap.~4, Lemma~1.6]
\label{thm:bound:above}
Let $f: \PP^n \dasharrow \PP^n$ be rational map 
of degree $d$ defined over $\overline{\QQ}$. 
Then there exists a constant $c_2$ such that 
\[
h_{nv}(f(x)) \leq d \; h_{nv}(x) + c_2
\]
for all $x \in \PP^n(\overline{\QQ})\setminus I_f$. 
\end{Theorem}

\medskip
\setcounter{equation}{0}
\section{Geometric properties of regular polynomial automorphisms}
\label{sec:geometry:general}
In this section, we show \eqref{eqn:Henon:blowups} 
for regular polynomial automorphisms of $\Aff^2$. 
First we recall the definition of 
regular polynomial automorphisms of $\Aff^2$. 
Consider a polynomial automorphism of degree $d \geq 2$
of the form 
\[
f \begin{pmatrix} x \\ y \end{pmatrix} 
= \begin{pmatrix} p(x,y)  \\ q(x,y) \end{pmatrix},  
\]
where $p(x,y)$ and $q(x,y)$ are polynomials in two variables, and 
$d$ is the maximum of $\deg p$ and $\deg q$. 
Let $\overline{f}: \PP^2 \dasharrow \PP^2$ be the extension of 
$f$ given in homogeneous coordinates as
\[
\overline{f} \begin{bmatrix} X \\Y \\Z \end{bmatrix} 
=  \begin{bmatrix}  
Z^d p(X/Z, Y/Z) \\ Z^d q(X/Z, Y/Z) \\Z^d \end{bmatrix}.  
\]
Let $H$ denote the line at infinity. 
Then $\overline{f}$ has a unique point 
of indeterminacy on $H$, denoted by $\bold{p}$.  
Let $f^{-1}: \Aff^2 \to \Aff^2$ be the inverse of $f$, 
and $\overline{f^{-1}}: \PP^2 \dasharrow \PP^2$ be its 
extension. Then $\overline{f^{-1}}$ has a unique point 
of indeterminacy on $H$, denoted by $\bold{q}$. 
A polynomial automorphism of $\Aff^2$ is said to be 
{\em regular} if $\bold{p} \neq \bold{q}$. 

By elimination of indeterminacy, by successively blowing up points 
starting from $\bold{p} \in \PP^2$, 
we obtain a projective surface $W$ and a composite of 
blow-ups $\pi_W: W \to \PP^2$ such that 
$\overline{f} \circ \pi_W: W \dasharrow \PP^2$  becomes a morphism. 
We take $W$ so that the number of blow-ups needed 
for elimination of indeterminacy is minimal. 
Noting that $\pi_W$ induces an isomorphism 
$\pi_W^{-1}(\PP^2\setminus\{\bold{p}\}) \to \PP^2\setminus\{\bold{p}\}$, 
we take $\bold{q'} \in W$ with $\pi_W(\bold{q'}) = \bold{q}$. 
In a parallel way as for $\bold{p}$, 
$\overline{f^{-1}}\circ\pi_W: W \dasharrow \PP^2$ 
becomes a morphism after a finite number of blow-ups 
starting at $\bold{q'}$.  

To summarize, there is a projective surface $V$ 
obtained by successive blow-ups of $\PP^2$ at $\bold{p}$ 
and then successive blow-ups at $\bold{q}$ 
in a parallel way as for $\bold{p}$ such that, 
if $\pi: V \to \PP^2$ denotes the morphism of blow-ups,
$\overline{f}\circ\pi$ extends to 
a morphism $\varphi: V \to \PP^2$ and 
$\overline{f^{-1}}\circ\pi$ extends to 
a morphism $\psi: V \to \PP^2$.
As for $W$, we take $V$ so that the number of blow-ups needed 
for elimination of indeterminacy is minimal. 
\begin{equation}
\label{eqn:figure:of:V}
\xymatrix{
       && V \ar[dll]_{\psi} \ar[d]^{\pi} \ar[drr]^{\varphi} && \\
\PP^2   && \PP^2 
\ar@{-->}[ll]_{\overline{f^{-1}}} \ar@{-->}[rr]^{\overline{f}} 
&& \PP^2}
\end{equation}

Before stating the next theorem, we fix some notation 
and terminology. Let $\rho: Y \to X$ be a morphism 
of smooth projective surfaces. 
For an irreducible curve $C$ on $Y$, 
its push-forward is defined by 
\[
\rho_*(C) := 
\begin{cases}
\deg(\rest{\rho}{C}: C \to f(C)) \ f(C) & \text{(if $f(C)$ is a curve)}, \\
0 & \text{(if $f(C)$ is a point)}. 
\end{cases}
\]
This extends linearly to a homomorphism 
$\rho_*$ from divisors on $Y$ to divisors on $X$. 
For two divisors $Z_1, Z_2$, we write $Z_1 \geq Z_2$ if 
$Z_1 - Z_2$ is effective. 

\begin{Theorem}
\label{thm:regular:estimate:pre}
Let $f: \Aff^2 \to \Aff^2$ be a regular 
polynomial automorphism of degree $d \geq 2$. 
Let $H$ denote the line at infinity. 
Let $V$ be as in \eqref{eqn:figure:of:V}. 
Then, as a $\QQ$-divisor on $V$,  
\[
D := \varphi^*H + \psi^*H - \left(d + \frac{1}{d}\right) \pi^*H  
\]
is effective. 
\end{Theorem}

\Proof
The proof we present here, which simplifies the proof we gave 
in the initial draft, is due to Noboru Nakayama. 

As above, let $\pi_W: W \to \PP^2$ be a composite of blow-ups of $\PP^2$ 
starting at $\bold{p}$ such that 
$\varphi_W :=\overline{f}\circ \pi_W: W \dasharrow \PP^2$ 
becomes a morphism. 
Let $H_W$ be the proper transform of $H$ by $\pi_W$, 
and $E_W$ the exceptional curve on $W$ given by 
the last blow-up of $\pi_W$. Since $\varphi_W$ is a morphism 
and $W$ is taken so that the number of blow-ups is minimal, 
we see that $\varphi_W$ sends 
$E_W$ to $H$ isomorphically. 

We consider $\pi_W^* H$ and $\varphi_W^* H$. We write 
$\pi_W^* H = a H_W + b E_W + M_W$ and  
$\varphi_W^* H = a' H_W + b' E_W + I_W$, 
where $a, b, a', b'$ are non-negative integers, and $M_W, I_W$ 
are effective divisors on $W$ with $\Supp(E_W) \not\subseteq \Supp(M_W), 
\Supp(E_W) \not\subseteq \Supp(I_W)$ 
such that $M_W, I_W$ are contracted to $\bold{p}$ 
by $\pi_W$. 

We determine $a, b, a', b'$. 
Since $\pi_W$ is a birational morphism, $\pi_{W*} \pi_W^* H = H$. 
It follows that $a = 1$. Similarly, $\varphi_{W*} \varphi_W^* H = H$ 
yields $b' =1$. On the other hand, let $[H]$ denotes 
the cohomology class of $H$ in $H^2(\PP^2, \ZZ)$. 
Since the degree of $f: \Aff^2 \to \Aff^2$ is $d$, 
we get $\varphi_{W*} \pi_W^* [H] = d [H] \in H^2(\PP^2, \ZZ)$. 
It follows that $\varphi_{W*} \pi_W^* H = d H$ and 
$b = d$. Since the degree of $f^{-1}: \Aff^2 \to \Aff^2$ 
is also $d$, we get 
$\pi_{W*} \varphi_W^* H = d H$ and $a' = d$. Putting together, we have
\begin{align*}
\pi_W^* H & = H_W + d E_W + M_W, \\ 
\varphi_W^* H & = d H_W + E_W + I_W.  
\end{align*}

Since the effective divisor $\pi_{W}^* H$ is nef, 
Lemma~\ref{lemma:N} below yields that 
\[
\varphi_{W}^* (d H) = \varphi_{W}^* (\varphi_{W*} \pi_W^* H) 
= (\varphi^*_{W} \varphi_{W*}) \pi_W^* H 
\geq \pi_W^* H. 
\]
We thus get
\begin{equation}
\label{eqn:N:1} 
d I_W \geq M_W.
\end{equation} 

In a parallel way as for $\bold{p}$, 
let $\pi_U: U \to \PP^2$ be a composite of blow-ups of $\PP^2$ 
starting at $\bold{q}$ such that 
$\psi_U :=\overline{f^{-1}} \circ \pi_U: U \dasharrow \PP^2$ 
becomes a morphism. 
Let $H_U$ be the proper transform of $H$ by $\pi_U$, 
and $F_U$ the exceptional curve on $U$ given by 
the last blow-up of $\pi_U$. The morphism $\psi_U$ sends 
$F_U$ to $H$ isomorphically. 
In a parallel way, we get 
\begin{align}
\notag \pi_U^* H & = H_U + d F_U + N_U, \\ 
\notag \psi_U^* H & = d H_U + F_U + J_U, \\
\label{eqn:N:2}
d J_U & \geq N_U,
\end{align}
where $N_U, J_U$ are effective divisors on $U$ 
with $\Supp(F_U) \not\subseteq \Supp(N_U), 
\Supp(F_U) \not\subseteq \Supp(J_U)$ 
such that $N_U, J_U$ are contracted to $\bold{q}$ 
by $\pi_U$. 

By the construction of $V$, there are birational morphisms 
$\alpha: V \to W$ and $\beta: V \to U$ such that the following 
diagram is commutative. 
\[
\xymatrix{
&& V \ar[dl]_{\beta} \ar[dd]^{\pi} \ar[dr]^{\alpha} 
\ar@/^2pc/[ddrr]^{\varphi} \ar@/_2pc/[ddll]_{\psi}& \\
&U \ar[dl]_{\psi_U} \ar[dr]^{\pi_U} & 
& W \ar[dl]_{\pi_W} \ar[dr]^{\varphi_W}\\
\PP^2 && \PP^2 \ar@{-->}[ll]_{\overline{f^{-1}}} 
\ar@{-->}[rr]^{\overline{f}} && \PP^2
}
\]

\smallskip
Let $H^\#$ on $V$ be the proper transform of $H$ by $\pi$. 
Let $E, M, I$ on $V$ be the proper transforms of $E_W, M_W, I_W$ 
by $\alpha$, respectively. Let $F, N, J$ be the proper transforms 
of $F_U, N_U, J_U$ by $\beta$, respectively. Then the following 
equalities hold: 
\begin{align}
\label{eqn:N:3}
\pi^* H & = H^\# + d E + d F + M + N, \\
\label{eqn:N:4}
\varphi^* H & = d (H^\# + d F + N) + E + I,  \\
\label{eqn:N:5}
\psi^* H & = d (H^\# + d E + M) + F + J.  
\end{align} 
By \eqref{eqn:N:3}--\eqref{eqn:N:5}, we get 
\begin{align*}
D & = \varphi^*H + \psi^*H - \left(d + \frac{1}{d}\right) \pi^*H \\
& = \left(d - \frac{1}{d}\right) H^\# 
- \frac{1}{d} M + I - \frac{1}{d} N + J. \\
\end{align*}
Since $d I \geq M$ and $d J \geq N$ 
by \eqref{eqn:N:1} and \eqref{eqn:N:2}, 
we see that $D$ is effective.
\QED

\begin{Lemma}
\label{lemma:N}
Let $\rho: Y \to X$ be a birational morphism of smooth 
projective surfaces. Let $Z$ be an effective divisor on $Y$. 
If $Z$ is nef, then $\rho^* \rho_* Z \geq Z$. 
\end{Lemma}

\Proof
First we treat a case when $\rho$ is the blow-up of $X$ 
at a point $x \in X$. Let $E$ denote the exceptional curve on 
$Y$. We write $Z = a_1 C_1 + \cdots + a_k C_k + b E$, 
where $C_1, \cdots, C_k, E$ are distinct irreducible and reduced 
curves, and $a_1, \cdots, a_k, b$ are non-negative integers.  
Then $\rho_* Z = a_1 \rho(C_1) + \cdots + a_k \rho(C_k)$. 
Hence $\rho^* \rho_* Z = a_1 (C_1 + m_1 E) 
+ \cdots + a_k (C_k + m_k E)$, where $m_i$ is the 
multiplicity of the curve $\rho(C_i)$ at $x$. 
Note that $m_i = C_i \cdot E$.  

Since $Z$ is nef, we get
\begin{align*}
Z \cdot E 
& = a_1 (C_1 \cdot E) + \cdots + a_1 (C_k \cdot E) + b (E \cdot E) \\ 
& = a_1 m_1 + \cdots + a_k m_k - b \geq 0. 
\end{align*}
Hence $a_1 m_1 + \cdots a_k m_k \geq b$ 
and we get $\rho^* \rho_* Z \geq Z$. 

In general, we decompose $\rho$ into a composite of blow-ups: 
$\rho = \rho_l \circ \cdots \circ \rho_2 \circ \rho_1$, where each 
$\rho_i$ is a blow-up at a point. Put 
$\rho' := \rho_l \circ \cdots \circ \rho_2$, and $Z' := \rho_{1*}Z$. 
Since the projection formula yields 
$(\rho_{1*}Z) \cdot C = Z \cdot (\rho_{1}^* C)$ for any curve, 
we see that $Z'$ is nef. Then, by induction, 
$\rho^{\prime *} \rho^{\prime}_* Z' \geq Z'$. 
Pulling back by $\rho_1$, we get $\rho_1^* (\rho^{\prime *} 
\rho^{\prime}_* Z') \geq \rho_1^* Z'$. Thus
\[
\rho^* \rho_* Z 
= \rho_1^* \rho^{\prime *} \rho^{\prime}_* (\rho_{1*} Z)
\geq \rho_1^* (\rho_{1*} Z) \geq Z.  
\]
\QED

Now we prove \eqref{eqn:Henon:blowups}. 

\begin{Theorem}
\label{thm:regular:estimate}
Let $f: \Aff^2 \to \Aff^2$ be a regular polynomial automorphism 
of degree $d \geq 2$ defined over a number field $K$. 
Then, there exists a constant $c$ such that 
\[
h_{nv}(f(x)) + h_{nv}(f^{-1}(x)) \geq 
\left(d + \frac{1}{d}\right) h_{nv}(x) -c
\]
for all $x\in \Aff^2(\overline{K})$. 
\end{Theorem}

\Proof
We can prove Theorem~\ref{thm:regular:estimate} 
as in \cite[Theorem~3.1]{SiHenon}. 
We take $x \in \Aff^2(\overline{K})$.  
Since $\pi: V\to \PP^2$ gives an isomorphism 
$\rest{\pi}{\pi^{-1}(\Aff^2)}: \pi^{-1}(\Aff^2) \to \Aff^2$, 
there is a unique point $\widetilde{x} \in V$ with 
$\pi(\widetilde{x}) = x$.  
By Theorem~\ref{thm:regular:estimate:pre}, we have 
\[
h_{V, \OO_V(\varphi^*H)}(\widetilde{x}) + 
h_{V, \OO_V(\psi^*H)}(\widetilde{x}) 
=
\left(d + \frac{1}{d}\right)h_{V, \OO_V(\pi^*H)}(\widetilde{x})
+ h_{V, \OO_V(D)}(\widetilde{x}) + O(1).
\]
It follows from Theorem~\ref{thm:height:machine}(1) that  
\[
h_{V, \OO_V(\varphi^*H)}(\widetilde{x}) 
= h_{\PP^2, \OO_V(H)}(\varphi(\widetilde{x})) + O(1)
= h_{\PP^2, \OO_V(H)}(f(x)) + O(1).
\] 
We similarly have  
\begin{align*}
h_{V, \OO_V(\psi^*H)}(\widetilde{x}) 
& = h_{\PP^2, \OO_V(H)}(f^{-1}(x)) + O(1),  \\
h_{V, \OO_V(\pi^*H)}(\widetilde{x})
& = h_{\PP^2, \OO_V(H)}(x) + O(1).
\end{align*}  
On the other hand, since $\pi(\Supp(D)) \subseteq \Supp(H)$,  
we have $\widetilde{x} \not\in \Supp(D)$.  
Since $D$ is effective by Theorem~\ref{thm:regular:estimate:pre}, 
it follows from Theorem~\ref{thm:height:machine}(2) that
there is a constant $c_2$ independent of $\widetilde{x}$ 
such that $h_{V, \OO_V(D)}(\widetilde{x}) \geq c_2$. 
Hence we get the assertion. 
\QED

\medskip
\setcounter{equation}{0} 
\section{H\'enon maps, conjugacy classes of polynomial automorphisms, 
and dynamical degrees}
\label{sec:geometry}
In this section, we review H\'enon maps, Friedland--Milnor's theorem 
on the conjugacy classes of polynomial automorphisms, and some properties 
of dynamical degrees of polynomial automorphisms, which will be used 
in \S\ref{sec:canonical:height:function}. 
We also give explicit forms of $\varphi^*H$, $\psi^*H$ and 
$\pi^*H$ in Theorem~\ref{thm:regular:estimate:pre} for H\'enon maps. 

A H\'enon map is a polynomial automorphism of the form
\begin{equation}
\label{eqn:Henon:maps}
f \begin{pmatrix} x \\ y \end{pmatrix} 
= \begin{pmatrix} p(x) - a y  \\ x \end{pmatrix},  
\end{equation}
where $a \neq 0$ and $p$ is a polynomial of degree $d \geq 2$. 
Let $\overline{f}: \PP^2 \dasharrow \PP^2$ 
(resp. $\overline{f^{-1}}: \PP^2 \dasharrow \PP^2$) 
be the birational extension of $f$ (resp. $f^{-1}$).  
Then $\overline{f}$ has the unique point of indeterminacy 
$\bold{p} = {}^t [0, 1, 0]$, and $\overline{f^{-1}}$ has the unique point 
of indeterminacy $\bold{q} = {}^t [1, 0, 0]$. 
In particular, H\'enon maps are examples of regular polynomial 
automorphisms. 

\medskip
We recall Friedland--Milnor's theorem \cite[\S2]{FM}, 
which is based on Jung's theorem \cite{Ju}. 
Let
\begin{equation}
\label{eqn:elementary}
E = \left\{\left.
f: \Aff^2 \to \Aff^2, 
\begin{pmatrix} x \\ y \end{pmatrix} 
\mapsto \begin{pmatrix} a x + P(y) \\ b y + c  \end{pmatrix}  
\;\right\vert\;
\begin{gathered}
a, b \in \overline{\QQ}^{\times}, c \in \overline{\QQ} \\
P(y) \in \overline{\QQ}[Y]
\end{gathered}
\right\}
\end{equation}
be the group of triangular automorphisms (also called 
de Jonqu\`eres automorphisms). 

\begin{Theorem}[\cite{FM}, \S2]
\label{thm:FM}
Let $f: \Aff^2 \to \Aff^2$ be a polynomial automorphism over $\overline{\QQ}$.
Then there is a polynomial automorphism $\gamma: \Aff^2 \to \Aff^2$ 
over $\overline{\QQ}$ such that $g := \gamma^{-1}\circ f \circ \gamma$ 
is one of the following types\textup{:}
\begin{enumerate}
\item[(i)]
$g$ is a triangular automorphism\textup{;}
\item[(ii)]
$g$ is a composite of H\'enon maps.
\end{enumerate}
\end{Theorem}

Note that Friedland--Milnor proved the theorem over $\CC$, but 
the theorem holds over $\overline{\QQ}$ by the specialization 
argument in \cite[Lemme~2]{De}. 

A polynomial automorphism $f$ is said to be {\em triangularizable} 
if it is conjugate to a triangular automorphism. 

\medskip
Here we recall dynamical degrees of polynomial automorphisms 
$f: \Aff^2 \to \Aff^2$. The dynamical degree of $f$ is defined by 
\[
\delta(f) := \lim_{n\to +\infty} \left(
\deg f^n \right)^{\frac{1}{n}} 
\]
(cf. \cite[D\'efinition~1.4.7]{Sib}). 
Suppose $g = \gamma^{-1}\circ f \circ \gamma$ is conjugate to $f$. 
Then, since $g^n = \gamma^{-1}\circ f \circ \gamma$, we have 
$\deg f^n - 2 \deg\gamma \leq \deg g^n \leq \deg f^n 
+ 2 \deg\gamma$. It follows that $\delta(f) = \delta(g)$. 
Thus dynamical degrees depend only on conjugacy classes 
of polynomial automorphisms. 

For polynomial automorphisms $g_1, g_2: \Aff^2 \to \Aff^2$ with degree 
$\deg g_1, \deg g_2 \geq 2$ and their extensions 
$\overline{g_1}, \overline{g_2}: \PP^2 \dasharrow \PP^2$, one has 
\begin{equation}
\label{eqn:criterion:of:being:regular}
\deg(g_1 \circ g_2) \leq (\deg g_1) (\deg g_2), 
\end{equation}
with equality if and only if the unique point $\bold{q}_{g_1}$ 
of indeterminacy of $\overline{g_1^{-1}}$ is different from the 
unique point $\bold{p}_{g_2}$ of indeterminacy of $\overline{g_2}$ 
(cf. \cite[Proposition~1.4.3]{Sib}).  
We remark that a composite $g$ of H\'enon maps 
is a regular polynomial automorphism, 
because the indeterminacy point of $\overline{g}$ is 
${}^t[0, 1, 0]$ while the indeterminacy point 
of $\overline{g^{-1}}$ is ${}^t [1, 0, 0]$.  

The following proposition is well-known. 

\begin{Proposition}
\label{prop:FM:supplimentary}
Let $f: \Aff^2 \to \Aff^2$ be a polynomial automorphism. 
Let $d$ be the degree of $f$ and $\delta$ 
the dynamical degree of $f$. 
\begin{enumerate}
\item[(1)]
$\delta$ is an integer with $1 \leq \delta \leq d$.  
\item[(2)]
$\delta =1$ if and only if $f$ is triangularizable.  
\item[(3)]
Suppose $d \geq 2$. Then 
$\delta = d$ if and only if $f$ is a regular 
polynomial automorphism. 
\end{enumerate}
\end{Proposition}

\Proof
We rely on the results of Furter \cite{Fu} to give a quick proof. 
We put $\tau = \frac{\deg(f^2)}{\deg f}$. Then 
Furter showed that either (i) $\tau \leq 1$ or (ii) $\tau$ is 
an integer greater than or equal to $2$. Moreover, (i) occurs 
if and only if $f$ is triangularizable 
(\cite[Propositon~5]{Fu}). 
In the case (ii), one has $\deg f^n = \tau^n \cdot \deg f$ 
(\cite[Propositon~4]{Fu}). 

(1) In the case (i), $f$ is triangular, and then  
its definition \eqref{eqn:elementary} yields that $\deg f^n \leq \deg f$, 
whence $\delta(f)=1$.  
In the case (ii), the dynamical degree of $f$ is equal to an 
integer $\tau \geq 2$. 

(2) It follows from the above proof of (1). 

(3) Since $d$ is assumed to be $\geq 2$, 
\eqref{eqn:criterion:of:being:regular} shows that 
$f$ is a regular polynomial automorphism if and only 
if $\tau = \deg f \ (\geq 2)$. 
Since $\tau = \delta(f)$ if $\tau \geq 2$, we get the assertion. 
\QED

\medskip
Since H\'enon maps are basic objects in the dynamics of 
polynomial automorphisms of $\Aff^2$ (cf. Theorem~\ref{thm:FM}), 
it would be worth giving explicit forms of $\varphi^*H$, $\psi^*H$ and 
$\pi^*H$ in Theorem~\ref{thm:regular:estimate:pre} for H\'enon maps
of degree $d \geq 2$, 
as Silverman \cite{SiHenon} did for quadratic H\'enon maps. 
In particular, this gives a different proof of 
Theorem~\ref{thm:regular:estimate:pre} in case of H\'enon maps. 

For this, we need an explicit description of blow-ups at
(infinitely near) points on $\PP^2$ that resolve the point of
indeterminacy of a H\'enon map $\overline{f}$. The case $\deg g = 2$ 
was carried out by Silverman \cite[\S2]{SiHenon}, and 
the general case by Hubbard--Papadopol--Veselov 
\cite[\S2]{HPV} in their compactification of H\'enon maps 
in $\CC^2$ as dynamical systems. 
Let us put together their results in the following theorem. 
(Note that, for the next theorem, the field of definition of 
$f$ can be any field, and $p(x)$ need not be monic.)

\begin{Theorem}[\cite{HPV}, \S2]
\label{thm:Henon:blowups}
\begin{enumerate}
\item[(1)]
Let $f$ be a H\'enon map in \eqref{eqn:Henon:maps}, 
and $\overline{f}: \PP^2 \dasharrow \PP^2$ its birational extension. 
Then $\overline{f}$ becomes 
well-defined after a sequence of $2d-1$ blow-ups. 
Explicitly, blow-ups are described as follows\textup{:}
\begin{enumerate}
\item[(i)]
First blow-up at ${\bold p}$\textup{;} 
\item[(ii)]
Next blow up at the unique point of indeterminacy, 
which is given by the intersection of the exceptional 
divisor and the proper transform of $H$\textup{;} 
\item[(iii)]
For the next $d-2$ times after \textup{(ii)}, 
blow-up at the unique point of indeterminacy, which 
is given by the intersection of the last exceptional divisor 
and the proper transform of the first exceptional 
divisor\textup{;}  
\item[(iv)]
For the next $d-1$ times after \textup{(iii)}, 
blow-up at the unique point of indeterminacy, 
which lies on the last exceptional divisor but 
not on the proper transform of the other exceptional 
divisors.  
\end{enumerate}
\item[(2)] 
Let $\overline{f_{2d-1}}: W \to \PP^2$ be the extension 
of the H\'enon map after the sequence of $2d-1$ blow-ups. 
Let $E_i^{'}$ denote the proper transform 
of $i$-th exceptional divisor \textup{(}$i=1, \cdots, 2d-1$\textup{)}. 
Then $\overline{f_{2d-1}}$ maps 
$E_i^{'}$ \textup{(}$i=1, \cdots, 2d-2$\textup{)} to ${\bold q}$, 
while $E_{2d-1}^{'}$ is mapped to $H$ by an isomorphism. 
\item[(3)] 
${E_1^{'}}^2 = -d$, ${E_i^{'}}^2 = -2$ 
\textup{(}$i=2, \cdots, 2d-2$\textup{)}, 
and ${E_{2d-1}^{'}}^2 = -1$. 
\end{enumerate}
\end{Theorem}

In particular, for H\'enon maps, $V$ in \eqref{eqn:figure:of:V} 
is the projective surface obtained by successive $2d-1$ 
blow-ups of $\PP^2$ at $\bold{p}$ as in 
Theorem~\ref{thm:Henon:blowups} and then 
successive $2d-1$ blow-ups at $\bold{q}$ 
in a parallel way as in Theorem~\ref{thm:Henon:blowups}. 

\begin{figure}
\begin{center}
\begin{picture}(320,240)
\put(0,20){\line(1,0){320}} 
\put(30,15){\line(2,1){50}}     %E_2
\put(80,30){\line(-2,1){40}}  %E_3
\put(40,40){\line(2,1){40}}    %E_4
\put(60,60){\circle*{1}}
\put(60,65){\circle*{1}}
\put(60,70){\circle*{1}}
\put(40,80){\line(2,1){40}}    %E_{d-1}
\put(70,85){\line(0,1){60}}   %E_{d}
\put(80,115){\line(-1,0){55}}  %E_1
\put(80,130){\line(-2,1){40}} %E_{d+1}
\put(60,150){\circle*{1}}
\put(60,155){\circle*{1}}
\put(60,160){\circle*{1}}
\put(80,170){\line(-2,1){40}} %E_{2d-3}
\put(40,180){\line(2,1){40}}   %E_{2d-2}
\put(80,190){\line(-2,1){40}}   %E_{2d-1}
%%%%%
\put(230,15){\line(2,1){50}}     %F_2
\put(280,30){\line(-2,1){40}}  %F_3
\put(240,40){\line(2,1){40}}    %F_4
\put(260,60){\circle*{1}}
\put(260,65){\circle*{1}}
\put(260,70){\circle*{1}}
\put(240,80){\line(2,1){40}}    %F_{d-1}
\put(270,85){\line(0,1){60}}   %F_{d}
\put(280,115){\line(-1,0){55}}  %F_1
\put(280,130){\line(-2,1){40}} %F_{d+1}
\put(260,150){\circle*{1}}
\put(260,155){\circle*{1}}
\put(260,160){\circle*{1}}
\put(280,170){\line(-2,1){40}} %F_{2d-3}
\put(240,180){\line(2,1){40}}   %F_{2d-2}
\put(280,190){\line(-2,1){40}}   %F_{2d-1}
%%%%%
\put(140,5){$H^{\#}$}
\put(25,5){$E_2$}
\put(85,25){$E_3$}
\put(25,35){$E_4$}
\put(25,75){$E_{d-1}$}
\put(65,75){$E_{d}$}
\put(10,110){$E_{1}$}
\put(85,125){$E_{d+1}$}
\put(85,165){$E_{2d-3}$}
\put(25,175){$E_{2d-2}$}
\put(85,185){$E_{2d-1}$}
\put(225,5){$F_2$}
\put(285,25){$F_3$}
\put(225,35){$F_4$}
\put(225,75){$F_{d-1}$}
\put(265,75){$F_{d}$}
\put(210,110){$F_{1}$}
\put(285,125){$F_{d+1}$}
\put(285,165){$F_{2d-3}$}
\put(225,175){$F_{2d-2}$}
\put(285,185){$F_{2d-1}$}
\end{picture}
\end{center}
\caption{The configuration after blow-ups.
The line $H^{\#}$ has the self-intersection number $-3$. 
The lines $E_1$ and $F_1$ have the self-intersection 
numbers $-d$. The lines $E_2, E_3, \cdots, E_{2d-2}$ 
and $F_2, F_3, \cdots, F_{2d-2}$ have 
the self-intersection numbers $-2$. The lines 
$E_{2d-1}$ and $F_{2d-1}$ have the self-intersection numbers $-1$.}
\label{figure}
\end{figure}
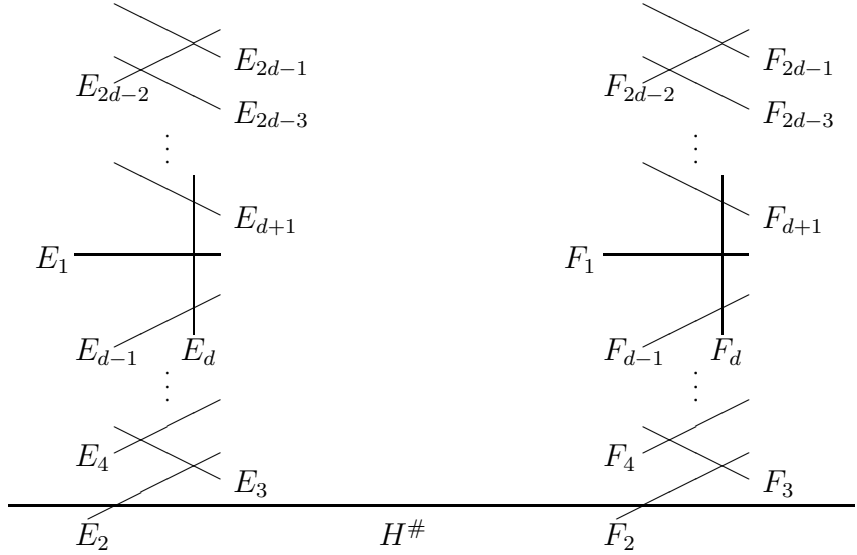

Let $E_i$ ($1\leq i\leq 2d-1$) be the proper transform 
of $i$-th exceptional divisor on $V$ on the side of $\bold{p}$, 
and $F_j$ ($1\leq j\leq 2d-1$) be the proper transform 
of $j$-th exceptional divisor on $V$ on the side of $\bold{q}$.  
Let $H^{\#}$ be the proper transform of $H$. 
The configuration of $H^{\#}$, $E_i$ and $F_j$ is 
illustrated in Figure~\ref{figure}. 

\begin{Proposition}
\label{prop:Henon:estimate}
Let $f: \Aff^2 \to \Aff^2$ be a H\'enon map of degree $d\geq 2$. 
Let the notation be as above. 
\begin{enumerate}
\item[(1)]
As divisors on $V$, we have  
{\allowdisplaybreaks
\begin{align*}
\pi^*H 
& = H^{\#} + \sum_{i=1}^{d} i E_i + \sum_{i=d+1}^{2d-1} d E_i 
+ \sum_{j=1}^{d} j F_j + \sum_{j=d+1}^{2d-1} d F_j,
\\
\varphi^*H 
& = d H^{\#} 
+ E_1 + \sum_{i=2}^{d} d E_i + \sum_{i=d+1}^{2d-1} (2d-i) E_i 
+ \sum_{j=1}^{d} j d F_j + \sum_{j=d+1}^{2d-1} d^2 F_j,\\
\psi^*H 
& = d H^{\#} 
+ \sum_{i=1}^{d} i d E_i + \sum_{i=d+1}^{2d-1} d^2 E_i
+ F_1 + \sum_{j=2}^{d} d F_j + \sum_{j=d+1}^{2d-1} (2d-j) F_j. 
\end{align*}}
\item[(2)]
The effective $\QQ$-divisor $D$ in 
Theorem~\ref{thm:regular:estimate:pre} is expressed as 
\begin{multline*}
D = \frac{d^2-1}{d} H^{\#} + 
\frac{d-1}{d} E_1 + \sum_{i=2}^d \frac{d^2 - i}{d} E_i + 
\sum_{i=d+1}^{2d-1} (2d-i-1) E_i \\
+ \frac{d-1}{d} F_1 + \sum_{j=2}^d \frac{d^2 - j}{d} F_j + 
\sum_{j=d+1}^{2d-1} (2d-j-1) F_j.
\end{multline*}
\end{enumerate}
\end{Proposition}

\Proof
We will show the expression for $\varphi^*H$. 
Since $\varphi$ maps $H^{\#}$, $E_i$ ($1 \leq i \leq 2d-2$) 
and $F_j$ ($1 \leq j \leq 2d-1$) to the point $\bold{q}$, 
we have 
\[
\varphi^*H \cdot H^{\#}=0, 
\qquad \varphi^*H \cdot E_i =0, 
\qquad \varphi^*H \cdot F_j =0 
\]   
for $1 \leq i \leq 2d-2$ and $1 \leq j \leq 2d-1$. 
Since $\varphi$ maps $E_{2d-1}$ to $H$ isomorphically, 
we have 
\[
\varphi^*H \cdot E_{2d-1} =1. 
\] 
Noting that the Picard group of $V$ is generated by 
$H^{\#}, E_i, F_j$ ($1 \leq i, j \leq 2d-1$), 
we set $\varphi^*H = a H^{\#} + \sum_{i=1}^{2d-1} b_i E_i 
+ \sum_{j=1}^{2d-1} c_j F_j$. 
From the above information and the information 
of the configuration after blow-ups (cf. Figure~\ref{figure}), 
we have the system of linear equations
\begin{equation*}
-3a + b_2 + c_2 =0, \qquad
\left\{ \begin{aligned} 
-d b_1 + b_d & = 0  \\
 a - 2 b_2 + b_3 & =0  \\
 b_{i-1} - 2 b_i + b_{i+1} & =0  \\
 b_1 + b_{d-1} - 2 b_{d} + b_{d+1}  & = 0  \\
 b_{2d-2} - b_{2d-1}  & =1, 
\end{aligned} \right. \qquad
\left\{ \begin{aligned}
 -d c_1 + c_d  & = 0  \\
 a - 2 c_2 + c_3  & =0  \\
 c_{j-1} - 2 c_j + c_{j+1} & =0 \\  
 c_1 + c_{d-1} - 2 c_{d} + c_{d+1}  & = 0  \\
 c_{2d-2} - c_{2d-1}  & =0, 
\end{aligned} \right.
\end{equation*}
where $i=3, \cdots, d-1, d+1, \cdots, 2d-2$ and 
$j=3, \cdots, d-1, d+1, \cdots, 2d-2$. 
By solving this system, we obtain the expression for $\varphi^*H$. 
Similarly we obtain the formula for $\psi^*H$. 
The formula for $\pi^*H$ follows from 
the construction of $V$. (We can also show this 
by using $\pi^*H \cdot H^{\#}=1$,   
$\pi^*H \cdot E_i =0$ and $\pi^*H \cdot F_j =0$ 
for all $i$ and $j$.)
The assertion (2) follows from (1). 
\QED

\begin{Remark}
Using classical results of Jung \cite{Ju} and van der Kulk \cite{Ku}, 
it is possible to explicitly compute $D$ for any regular polynomial 
automorphisms $f$ of degree $d \geq 2$, as in 
Proposition~\ref{prop:Henon:estimate} for H\'enon maps. In this case, 
coefficients of $D$ are expressed in terms of the 
polydegree $(d_1, \ldots, d_l)$ of $f$ (cf. \cite[\S3]{FM}). 
Note that, for H\'enon maps $f$ of degree $d \geq 2$, 
its polydegree is $(d)$, i.e., $l=1$ and $d_1 = d$.   
\end{Remark}

\medskip
\setcounter{equation}{0}
\section{Canonical height functions}
\label{sec:canonical:height:function}
In this section, we will prove Theorem~\ref{thm:main} and 
Corollary~\ref{cor:main} by showing Theorem~\ref{thm:general}. 
We first fix some notation and terminology. 
We refer to the survey \cite{Sib} for more details 
about the dynamics of polynomial automorphisms.  

Let $f: \Aff^n \to \Aff^n$ be a polynomial automorphism 
over a number field $K$. 
We use the notation $\overline{f}$ 
to denote the birational extension of $f$ to $\PP^n$. 
Let $f^{-1}: \Aff^n \to \Aff^n$ denote the inverse of $f$, 
and we use the notation $\overline{f^{-1}}$ to denote 
the birational extension of $f^{-1}$ to $\PP^n$. 
Note that the degree of $f$ and the degree of $f^{-1}$ 
may not be the same when $n \geq 3$ (cf. \cite[Chapitre~2]{Sib}). 

Let $S$ be a set and $T$ a subset of $S$. 
Two real-valued functions $\lambda$ and $\lambda'$ 
on $S$ are said to be equivalent on $T$ if there exist 
positive constants $a_1$, $a_2$ 
and constants $b_1$, $b_2$ such that 
$a_1 \lambda(x) + b_1 \leq \lambda'(x) \leq a_2 \lambda(x) + b_2$ 
for all $x \in T$. 
We use the notation 
$\lambda \gg\ll \lambda'$
to denote this equivalence. 
(Note that our notation $\gg\ll$ is different from 
that in \cite[Chap.~4, \S1]{La} where $b_1 = b_2 = 0$. )

\begin{Theorem}
\label{thm:general}
Let $f: \Aff^n \to \Aff^n$ be a polynomial automorphism over 
a number field $K$. Let $\gamma: \Aff^n \to \Aff^n$ be 
a polynomial automorphism over $K$, 
and we define the polynomial automorphism 
$g: \Aff^n \to \Aff^n$ by $g := \gamma^{-1} \circ f \circ \gamma$. 
Let $\delta$ and $\delta_-$ denote the degrees of $g$ and 
$g^{-1}$, respectively. 
We assume that $\delta \geq 2$ and 
that there exists a constant $c$ such that 
\begin{equation}
\label{eqn:general}
\frac{1}{\delta} h_{nv}(g(x)) + \frac{1}{\delta_-} h_{nv}(g^{-1}(x))
\geq \left(1 + \frac{1}{\delta \delta_-}\right) h_{nv}(x) - c 
\end{equation}
for all $x \in \Aff^n(\overline{K})$.  
Then there exists a function 
$\widehat{h}: \Aff^n(\overline{K}) \to \RR$
with the following properties\textup{:} 
\begin{enumerate}
\item[(i)] 
$h_{nv} \gg\ll \widehat{h}$ on $\Aff^n(\overline{K})$\textup{;}
\item[(ii)] 
$\frac{1}{\delta}\widehat{h} \circ f + 
\frac{1}{\delta_-}\widehat{h} \circ f^{-1}
= \left(1 + \frac{1}{\delta \delta_-} \right) \widehat{h}$. 
\end{enumerate}
Moreover, $\widehat{h}$ enjoys the following uniqueness property\textup{:}  
if $\widehat{h}'$ is another function satisfying (i) and (ii) such that 
$\widehat{h}' = \widehat{h} + O(1)$, then 
$\widehat{h}' = \widehat{h}$. 
Furthermore, 
$\widehat{h}(x) \geq 0$ for all $x \in \Aff^n(\overline{K})$,  
and $\widehat{h}(x) = 0$ if and only if $x$ is $f$-periodic. 
\end{Theorem}

\medskip
{\sl Proof of Theorem~\ref{thm:main} and Corollary~\ref{cor:main}. }\quad
Admitting Theorem~\ref{thm:general}, we will prove 
Theorem~\ref{thm:main} and Corollary~\ref{cor:main}.
We may replace $K$ by a finite extension field. 
Since the dynamical degree $\delta$ is greater than or equal to $2$,  
Theorem~\ref{thm:FM} and Proposition~\ref{prop:FM:supplimentary} 
yield that there is a polynomial automorphism 
$\gamma$ such that $g := \gamma\circ f\circ \gamma^{-1}$ is 
a composite of H\'enon maps. Since a composite of H\'enon maps 
is a regular polynomial automorphism 
(cf. lines before Proposition~\ref{prop:FM:supplimentary}), 
it follows from Theorem~\ref{thm:regular:estimate} that 
$g$ satisfies \eqref{eqn:general}.
Then, noting that the dynamical degrees of $f$ and $g$ are the same, 
Theorem~\ref{thm:main} and Corollary~\ref{cor:main} 
follows from Theorem~\ref{thm:general}.
\QED

\medskip
{\sl Proof of Theorem~\ref{thm:general}. }\quad

{\bf Step 1.}\quad
We show the existence of a function 
$\widehat{h}_g: \Aff^n(\overline{K}) \to \RR$
with the following properties: 
\begin{enumerate}
\item[(iii)] 
$h_{nv} \gg\ll \widehat{h}_g$ on $\Aff^n(\overline{K})$\textup{;}
\item[(iv)] 
$\frac{1}{\delta}\widehat{h}_g \circ g + 
\frac{1}{\delta_-}\widehat{h}_g \circ g^{-1}
= \left(1 + \frac{1}{\delta \delta_-} \right) \widehat{h}_g$. 
\end{enumerate}

For $x \in \Aff^n(\overline{K})$, we define 
\[
\widehat{h}_{g}^{+}(x) 
= \limsup_{l \to  +\infty} \frac{1}{\delta^l} h_{nv} (g^l(x)), 
\qquad
\widehat{h}_{g}^{-}(x) = \limsup_{l \to  +\infty} 
\frac{1}{\delta_{-}^l} h_{nv} (g^{-l}(x)),  
\]
{\em a priori} in $\RR\cup \{\infty\}$, but we will show in the next claim 
that this value is finite. We define 
\[
\widehat{h}_{g}(x) = \widehat{h}_{g}^{+}(x) + \widehat{h}_{g}^{-}(x). 
\]
Note that 
this definition of $\widehat{h}^{\pm}_{g}$ has some similarity 
to the definition of Green currents on $\Aff^n(\CC)$ 
associated with $g$ (cf. \cite[D\'efinition~2.2.5]{Sib}),  
and to  Silverman's definition of canonical heights 
on certain K3 surfaces \cite[\S3]{SiK3}. 
Let us show $\widehat{h}_{g}$ satisfies the properties 
(iii) and (iv). 

\begin{Claim}
\label{claim:general:1}
There exist constants $c^{\pm}$ such that 
$\widehat{h}_{g}^{\pm}(x) \leq h_{nv}(x) + c^{\pm}$
for all $x \in \Aff^n(\overline{K})$. 
\end{Claim}

\Proof
By Theorem~\ref{thm:bound:above}, there exists a constant 
$c_2$ such that $\frac{1}{\delta} h_{nv}(g(x)) \leq h_{nv}(x) + \frac{c_2}{\delta}$ 
for all $x \in \Aff^n(\overline{K})$. 
We show 
\[
\frac{1}{\delta^l} h_{nv}(g^l(x)) 
\leq h_{nv}(x) + \left(\sum_{i=1}^l \frac{1}{\delta^i}\right) c_2
\]
by the induction on $l$. Indeed, since 
$\frac{1}{\delta} h_{nv}(g^{l+1}(x)) \leq h_{nv}(g^{l}(x)) + \frac{c_2}{\delta}$, 
we have 
\[
\frac{1}{\delta^{l+1}} h_{nv}(g^{l+1}(x)) 
\leq \frac{1}{\delta^{l}} h_{nv}(g^{l}(x)) + \frac{c_2}{\delta^{l+1}}
\leq h_{nv}(x) + \left(\sum_{i=1}^{l+1} \frac{1}{\delta^i}\right) c_2. 
\] 
By putting $c^{+} = c_2 \sum_{i=1}^{ +\infty} \frac{1}{\delta^i} 
= \frac{c_2}{\delta-1}$, we obtain 
$\widehat{h}_{g}^{+}(x) = \limsup_{l \to  +\infty} 
\frac{1}{\delta^l} h_{nv}(g^l(x)) \leq h_{nv}(x) + c^{+}$. 
The estimate for $\widehat{h}_{g}^{-}$ is shown similarly. 
(Note that it follows from $\delta\geq 2$ that $\delta_{-}\geq 2$.)
\QED

\begin{Claim}
\label{claim:general:2}
We have 
\[
\widehat{h}_{g}(x) \geq 
h_{nv}(x) - \frac{\delta \delta_-}{(\delta-1)(\delta_{-}- 1)}c 
\]
for all $x \in \Aff^n(\overline{K})$, 
where $c$ is the constant given in \eqref{eqn:general}. 
\end{Claim}

\Proof
We set $h' = h_{nv} - \frac{\delta \delta_-}{(\delta-1)(\delta_{-}- 1)}c$. 
Then we have for all $x \in \Aff^n(\overline{K})$ 
\begin{equation}
\label{eqn:general:2}
\frac{1}{\delta} h'(g(x)) + \frac{1}{\delta_-} h'(g^{-1}(x))
\geq \left(1 + \frac{1}{\delta \delta_-}\right) h'(x). 
\end{equation}
Then we have 
%\begin{gather*}
$\frac{1}{\delta^2} h'(g^2(x)) + \frac{1}{\delta \delta_-} h'(x)
\geq \left(1 + \frac{1}{\delta \delta_-}\right) \frac{1}{\delta}h'(g(x)) $
and 
$\frac{1}{\delta\delta_-} h'(x) + \frac{1}{\delta_-^2} h'(g^{-2}(x))
\geq \left(1 + \frac{1}{\delta \delta_-}\right) \frac{1}{\delta_-}h'(g^{-1}(x))$. 
%\end{gather*}
Adding these two inequalities and 
using \eqref{eqn:general:2} again, we obtain 
\[
\frac{1}{\delta^2} h'(g^2(x)) +  \frac{1}{\delta_-^2} h'(g^{-2}(x))
\geq  \left(1 + \frac{1}{(\delta \delta_-)^2}\right) h'(x).
\]
Inductively, we obtain 
\[
\frac{1}{\delta^{2^l}} h'(g^{2^l}(x)) +  \frac{1}{\delta_-^{2^l}} h'(g^{-{2^l}}(x))
\geq  \left(1 + \frac{1}{(\delta \delta_-)^{2^l}}\right) h'(x).
\]
(Though not necessary for the proof, one can also show 
$\frac{1}{\delta^{m}} h'(g^{m}(x)) +  \frac{1}{\delta_-^{m}} h'(g^{-{m}}(x))
\geq  \left(1 + \frac{1}{(\delta \delta_-)^{m}}\right) h'(x)$ 
for every $m\in\ZZ_{> 0}$.)
By letting $l \to  +\infty$, it follows that 
\begin{multline}
\label{eqn:general:3}
\limsup_{l \to  +\infty} \frac{1}{\delta^{2^l}} h'(g^{2^l}(x)) 
+ \limsup_{l \to  +\infty}  \frac{1}{\delta_-^{2^l}} h'(g^{-{2^l}}(x)) \\
\geq \limsup_{l \to  +\infty} \left(
\frac{1}{\delta^{2^l}} h'(g^{2^l}(x)) 
+  \frac{1}{\delta_-^{2^l}} h'(g^{-{2^l}}(x))
\right)
\geq h'(x).
\end{multline}
Since 
\begin{multline*}
\widehat{h}_{g}^+(x)
= 
\limsup_{m \to  +\infty} \frac{1}{\delta^{m}} h_{nv}(g^{m}(x)) \\
=
\limsup_{m \to  +\infty} \frac{1}{\delta^{m}} 
\left( h'(g^{m}(x)) +  \frac{\delta \delta_-}{(\delta-1)(\delta_{-}- 1)}c \right)
\geq 
\limsup_{l \to  +\infty} \frac{1}{\delta^{2^l}} h'(g^{2^l}(x)) 
\end{multline*} 
and similarly 
$\widehat{h}_{g}^-(x)
\geq 
\limsup_{l \to  +\infty} \frac{1}{\delta_-^{2^l}} h'(g^{-2^l}(x))$, 
the left-hand-side of \eqref{eqn:general:3} is less than 
or equal to $\widehat{h}_{g}(x)$, while the right-hand-side is 
$h_{nv}(x) - \frac{\delta \delta_-}{(\delta-1)(\delta_{-}- 1)}c$. 
Thus we get the desired inequality. 
\QED

The property (iii) follows from Claim~\ref{claim:general:1} 
and Claim~\ref{claim:general:2}. Indeed we have  
\begin{equation}
\label{eqn:general:4}
h_{nv}(x) - \frac{\delta \delta_-}{(\delta-1)(\delta_{-}- 1)}c
\leq \widehat{h}_{g}(x) \leq 
2 h_{nv}(x) + c^+ + c^-. 
\end{equation}
The property (iv) is checked by the following equations:
\begin{gather*}
\widehat{h}_{g}^+(f(x)) = \delta \widehat{h}_{g}^+(x), 
\quad \widehat{h}_{g}^+(f^{-1}(x)) 
= \frac{1}{\delta}\widehat{h}_{g}^+(x); \\
\widehat{h}_{g}^-(f(x)) = \frac{1}{\delta_-} \widehat{h}_{g}^-(x), 
\quad \widehat{h}_{g}^-(f^{-1}(x)) = \delta_- \widehat{h}_{g}^-(x). 
\end{gather*}
Thus $\widehat{h}_{g}: \Aff^n(\overline{K}) \to \RR$ 
satisfies the properties (iii) and (iv). 

\smallskip
{\bf Step 2.}\quad
We show the existence of a function 
$\widehat{h}_{\circ}: \Aff^n(\overline{K}) \to \RR$ with the properties 
(i) and (ii). We define $\widehat{h}_{\circ}$ by 
\[
\widehat{h}_{\circ}(x) := \widehat{h}_{g}(\gamma^{-1}(x))
\]
for all $x \in \Aff^n(\overline{K})$. 

By \eqref{eqn:general:4}, we have 
$\widehat{h}_{g}(\gamma^{-1}(x)) \leq 2 h_{nv}(\gamma^{-1}(x)) 
+ c^+ + c^-$. Theorem~\ref{thm:bound:above} yields that 
there is a constant $c_{\gamma^{-1}}$ such that 
$h_{nv}(\gamma^{-1}(x)) \leq (\deg\gamma^{-1})\ h_{nv}(x) 
+ c_{\gamma^{-1}}$ for all $x \in \Aff^n(\overline{K})$. Thus
\begin{equation}
\label{eqn:general:5}
\widehat{h}_{\circ}(x) \leq 2 (\deg\gamma^{-1})\ h_{nv}(x)   
+ (2 c_{\gamma^{-1}} + c^+ + c^-).  
\end{equation}
On the other hand, Theorem~\ref{thm:bound:above} 
yields that there is a constant $c_{\gamma}$ such that 
$h_{nv}(\gamma(x)) \leq (\deg\gamma)\ h_{nv}(x) 
+ c_{\gamma}$ for all $x \in \Aff^n(\overline{K})$. 
Hence 
\[
h_{nv}(\gamma^{-1}(x)) \geq (\deg\gamma)^{-1} h_{nv}(x) 
- (\deg\gamma)^{-1} c_{\gamma}. 
\]
Then by \eqref{eqn:general:4}, we get 
\begin{equation}
\label{eqn:general:6}
\widehat{h}_{\circ}(x) \geq 
(\deg\gamma)^{-1} h_{nv}(x) 
- (\deg\gamma)^{-1} c_{\gamma}
- \frac{\delta \delta_-}{(\delta-1)(\delta_{-}- 1)}c
\end{equation}
for all $x \in \Aff^n(\overline{K})$. 
Now the property (i) follows from 
\eqref{eqn:general:5} and \eqref{eqn:general:6}. 

The property (iv) follows from 
\begin{align*}
\widehat{h}_{\circ}(f(x)) + \widehat{h}_{\circ}(f^{-1}(x)) 
& = \widehat{h}_{g}(\gamma^{-1}(f(x))) 
+ \widehat{h}_{g}(\gamma^{-1}(f^{-1}(x))) \\
& = \widehat{h}_{g}(g(\gamma^{-1}(x))) 
+ \widehat{h}_{g}(g^{-1}(\gamma^{-1}(x))) \\
& = \left(1 + \frac{1}{\delta \delta_-} \right) 
\widehat{h}_g(\gamma^{-1}(x)) 
=  \left(1 + \frac{1}{\delta \delta_-} \right) 
\widehat{h}_{\circ}(x), 
\end{align*}
where we used (iv) in the third equality. 

\smallskip
{\bf Step 3.}\quad
We will show uniqueness property of $\widehat{h}$. 
In what follows, let $\widehat{h}$ denote a function 
with the properties (i) and (ii), not necessarily 
being equal to $\widehat{h}_{\circ}$. 

Suppose $\widehat{h}'$ is another function with the properties 
(i) and (ii) such that $\lambda:= \widehat{h}' - \widehat{h}$ 
is bounded on $\Aff^n(\overline{K})$. Set 
$M : = \sup_{x \in \Aff^n(\overline{K})}\vert{\lambda}(x)\vert$.  
Then 
\begin{multline*}
\left(1 + \frac{1}{\delta \delta_-} \right) M 
=  \left(1 + \frac{1}{\delta \delta_-} \right) 
\sup_{x \in \Aff^n(\overline{K})} \left\vert\lambda(x)\right\vert \\
= \sup_{x \in \Aff^n(\overline{K})} 
\left\vert\frac{1}{\delta}\lambda (f(x)) +  
\frac{1}{\delta_-}\lambda (f^{-1}(x)) \right\vert
\leq 
\left(\frac{1}{\delta} + \frac{1}{\delta_-} \right) M. 
\end{multline*}
Since $1 + \frac{1}{\delta \delta_-} - \frac{1}{\delta} - \frac{1}{\delta_-} 
= \frac{(\delta-1)(\delta_- -1)}{\delta \delta_-} > 0$, we have $M=0$, 
hence $\widehat{h} = \widehat{h}'$.  

To show $\widehat{h} \geq 0$, we assume the contrary,  
so that there exists $x_0 \in \Aff^n(\overline{K})$ 
with $\widehat{h}(x_0) =: a < 0$. 
Then $\frac{1}{\delta}\widehat{h}(f(x_0)) + 
\frac{1}{\delta_-}\widehat{h}(f^{-1}(x_0))
= \left(1 + \frac{1}{\delta \delta_-} \right) \widehat{h}(x_0) 
= \left(1 + \frac{1}{\delta \delta_-} \right) a$. 
Thus we have 
\[
\widehat{h}(f(x_0)) \leq \frac{1+\delta \delta_-}{\delta + \delta_-} a
\quad\text{or}\quad 
\widehat{h}(f^{-1}(x_0)) \leq \frac{1+\delta \delta_-}{\delta + \delta_-} a.
\]
Since $\frac{1+\delta \delta_-}{\delta + \delta_-} > 1$, this shows that 
$\widehat{h}$ is not bounded from below. 
Since $h_{nv}$ is bounded from below and 
$h_{nv} \gg\ll \widehat{h}$, this is a contradiction. 

Finally we will show that $x \in \Aff^n(\overline{K})$ is 
$f$-periodic if and only if $\widehat{h}(x)=0$. 

Suppose $\widehat{h}(x_1)=0$. Then by (ii) 
and the non-negativity of $\widehat{h}$, we have 
$\widehat{h}(f(x_1))=0$ and $\widehat{h}(f^{-1}(x_1))=0$. 
Take an extension field $L$ of $K$ such 
that $x_1$ is defined over $L$. 
Since $\widehat{h} \gg\ll h_{nv}$, 
$\widehat{h}$ satisfies the Northcott finiteness property. 
Thus the set
\[
\{ f^l(x_1) \mid l \in \ZZ \} 
\quad\left(\subseteq 
\{ x \in \Aff^n(L) \mid \widehat{h}(x) = 0 
\}\right) 
\]
is finite. Hence $x_1$ is $f$-periodic.  

On the other hand, suppose $\widehat{h}(x_2) =:b > 0$. 
Then it follows from (ii) that  
\[
\widehat{h}(f(x_2)) \geq \frac{1+\delta \delta_-}{\delta + \delta_-} b
\quad\text{or}\quad 
\widehat{h}(f^{-1}(x_2)) \geq \frac{1+\delta \delta_-}{\delta + \delta_-} b.
\]
This shows that the set $\{f^l(x_2) \mid l \in \ZZ\}$ 
is not a set of bounded height. Thus $x_2$ cannot be $f$-periodic.  
\QED

In the remainder of this section, we would like to discuss 
the condition \eqref{eqn:general} in Theorem~\ref{thm:general}. 
The next proposition shows that 
the constant $(1 + \frac{1}{\delta \delta_-})$ 
in \eqref{eqn:general} is the largest number one can hope for. 

\begin{Proposition}
\label{prop:upper:bound:of:coeff}
Let $g: \Aff^n \to \Aff^n$ a polynomial automorphism of 
degree $\delta \geq 2$ over a number field $K$. Let $\delta_{-}$ denote the degree 
of $g^{-1}$. Let $a \in \RR$. Suppose there exists a constant $c$ 
such that  
\[
\frac{1}{\delta} h_{nv}(g(x)) + \frac{1}{\delta_-} h_{nv}(g^{-1}(x))
\geq a h_{nv}(x) - c 
\]
for all $x \in \Aff^n(\overline{K})$. 
Then $a \leq 1 + \frac{1}{\delta \delta_-}$. 
\end{Proposition}

\Proof
To lead a contradiction, we assume that $a > 1 + \frac{1}{\delta \delta_-}$. 
Noting $a > 1 + \frac{1}{\delta \delta_-} \geq \frac{1}{\delta} + \frac{1}{\delta_-}$, 
we set $c' := \left( a - \frac{1}{\delta} - \frac{1}{\delta_-}\right)^{-1} c$ 
and $h' := h_{nv} -c'$. Then $h'$ satisfies
\begin{equation}
\label{eqn:upper:bound:of:coeff}
\frac{1}{\delta} h'(g(x)) + \frac{1}{\delta_-} h'(g^{-1}(x))
\geq a h'(x)
\end{equation}
for all $x \in  \Aff^n(\overline{K})$. 
As in the proof of Claim~\ref{claim:general:2}, 
we get
\[
\frac{1}{\delta^2} h'(g^2(x)) + \frac{1}{\delta_-^2} h'(g^{-2}(x))
\geq \left(a^2 - \frac{2}{\delta \delta_-}\right) h'(x). 
\]
We set $a_1 = a^2 - \frac{2}{\delta \delta_-}$. 
Since $a_1 - 1 - \frac{1}{(\delta \delta_-)^2} 
= a^2 - \frac{2}{\delta \delta_-} - 1 - \frac{1}{(\delta \delta_-)^2} 
> (1 + \frac{1}{\delta \delta_-})^2 - \frac{2}{\delta \delta_-} - 1 - \frac{1}{(\delta \delta_-)^2} 
= 0$, we have $a_1 > 1 + \frac{1}{(\delta \delta_-)^2}$. 
Thus, if we define a sequence $\{a_l\}_{l=0}^{ +\infty}$ 
by $a_0 = a$ and $a_{l+1} = a_l^2 - \frac{2}{(\delta \delta_-)^{2^l}}$, 
then we get inductively
\[
\frac{1}{\delta^{2^l}} h'(g^{2^l}(x)) + \frac{1}{\delta_-^{2^l}} 
h'(g^{-{2^l}}(x)) \geq a_{l} h'(x). 
\]
On the other hand, it follows from Theorem~\ref{thm:bound:above} 
and the argument in Claim~\ref{claim:general:1} 
that there is a constant $c''$ independent of 
$l \in \ZZ$ such that for all $x \in \Aff^2(\overline{K})$, 
\[
2 h'(x) + c'' \geq 
\frac{1}{\delta^{2^l}} h'(g^{2^l}(x)) + \frac{1}{\delta_-^{2^l}} 
h'(g^{-{2^l}}(x)).  
\]
Thus $2 h' + c'' \geq a_l h'$. 
Since $h' = h_{nv} -c'$ and 
$\lim_{l \to  +\infty} a_l =  +\infty$ follows from  
Lemma~\ref{lemma:sequence}(1), this is a contradiction.  
\QED

\begin{Lemma}
\label{lemma:sequence}
Let $D \geq 4$. Let $\{a_l\}_{l=0}^{ +\infty}$ be a sequence  
defined by $a_0 = a$ and $a_{l+1} 
= a_l^2 - 2 D^{-2^l}$. 
\begin{enumerate}
\item[(1)] 
If $a > 1 + \frac{1}{D}$, then $\lim_{l \to  +\infty} a_l =  +\infty$. 
\item[(2)] 
If $a = 1 + \frac{1}{D}$, then $\lim_{l \to  +\infty} a_l = 1$.
\item[(3)] 
If $1 \leq a < 1 + \frac{1}{D}$, then $\lim_{l \to  +\infty} a_l = 0$.
\end{enumerate}
\end{Lemma}

\Proof
We show (1). Set $\varepsilon_l = a_l - 1 - D^{-2^l}$. 
In particular $\varepsilon_0 = a - 1 - D^{-1} > 0$. 
Since 
$\varepsilon_{l+1} 
=  a_{l+1} - 1 - D^{-2^{l+1}}
=  2 \varepsilon_l (1 + D^{-2^l}) + \varepsilon_l^2$, 
we get $\varepsilon_{l+1} > 2 \varepsilon_l 
> \cdots > 2^{l+1}\varepsilon_0$. Hence 
$\lim_{l\to +\infty} \varepsilon_{l} =  +\infty$ 
and thus $\lim_{l\to +\infty} a_{l} =  +\infty$

We show (2). In this case, 
we have $a_l = 1 + D^{-2^l}$. Thus $\lim_{l\to +\infty} a_{l} = 1$.

Finally we show (3). On one hand, we get by induction 
$a_l \geq 2 D^{-2^{l-1}}$ for $l\geq 1$, and in particular  
$a_l \geq 0$ for $l\geq 1$. 
On the other hand, we claim for sufficiently large $l$ 
that $a_l < 1$. Indeed, we assume the contrary and suppose
$a_l \geq 1$ for all $l$.
By induction, we get $a_l < 1 + D^{-2^l}$. 
We set $\lambda_l = 1 + D^{-2^l} - a_l$, 
and so $0 < \lambda_l \leq D^{-2^l}$. 
Then 
$a_{l+1} = a_l^2 - 2 D^{-2^{l}} = 
(1 + D^{-2^l} - \lambda_l)^2 - 2 D^{-2^{l}}
= 1 + D^{-2^{l+1}} - 2 \lambda_l (1 + D^{-2^l}) +\lambda_l^2$. 
Hence we get 
$\lambda_{l+1} = 2 \lambda_l (1 + D^{-2^l}) - 
\lambda_l^2 \geq 2 \lambda_l$,  
which says that $\lim_{l\to +\infty} \lambda_l =  +\infty$. 
This is a contradiction. Hence there is an $l_0$ 
with $a_{l_0} < 1$. Since $(0 \leq)\; a_{l_0+k} 
\leq a_{l_0}^{2^k}$, we get $\lim_{l\to +\infty} a_l = 0$. 
\QED

Let $a_{\sup}$ denote the supremum of $a \in \RR$ 
that satisfies the inequality in 
Proposition~\ref{prop:upper:bound:of:coeff}. 
It follows from Theorem~\ref{thm:regular:estimate} that,  
if $g$ is a regular polynomial automorphism of $\Aff^2$ 
of degree $\delta \geq 2$, then $\delta=\delta_-$ and 
$a_{\sup} = 1 + \frac{1}{\delta^2}$.
We remark that 
Marcello \cite[Th\'eor\`eme~3.1]{Ma} 
showed that, if $g$ is a regular polynomial automorphism 
of $\Aff^n$ (this means the set of indeterminacy 
$I_{\overline{g}}$ and  $I_{\overline{g^{-1}}}$ are disjoint, 
cf. \cite[D\'efinition~2.2.1]{Sib}), 
then $a_{\sup} \geq 1$. 
It would be interesting to know what polynomial automorphisms $g$ 
on $\Aff^n$ satisfy \eqref{eqn:general}.

\medskip
\setcounter{equation}{0}
\section{The number of points with bounded height in an $f$-orbit}
\label{sec:estimate}
In this section, we will prove Theorem~\ref{thm:estimate}. 
As in \S\ref{sec:canonical:height:function}
we will show Theorem~\ref{thm:estimate} 
in a more general setting. The arguments below are inspired by 
those of Silverman on certain K3 surfaces \cite[\S3]{SiK3}. 

{\em Throughout this section, 
let $f: \Aff^n \to \Aff^n$ be a polynomial automorphism 
of over a number field $K$ satisfying the conditions 
in Theorem~\ref{thm:general}. Let $\widehat{h}$ be a height 
function constructed in Theorem~\ref{thm:general}.}

We define functions $\widehat{h}^{\pm}: \Aff^n(\overline{K}) \to \RR$ 
to be 
\begin{align*}
\widehat{h}^+(x) &
= \frac{\delta \delta_-}{(\delta \delta_-)^2 - 1} 
\left(\delta_- \widehat{h}(f(x)) - \frac{1}{\delta_-}\widehat{h}(f^{-1}(x))\right), \\
\widehat{h}^-(x) & 
= \frac{\delta \delta_-}{(\delta \delta_-)^2 - 1} 
\left(\delta \widehat{h}(f^{-1}(x)) 
- \frac{1}{\delta}\widehat{h}(f(x))\right) 
\end{align*}
for $x \in \Aff^n(\overline{K})$. 
We remark that, in the notations of the proof 
of Theorem~\ref{thm:general}, if $\widehat{h} = \widehat{h}_g$, 
then $\widehat{h}^+ = \widehat{h}_g^+$ and 
$\widehat{h}^- = \widehat{h}_g^-$.

\begin{Lemma}
\label{lemma:h:+:-}
\begin{enumerate}
\item[(1)]
$\widehat{h} = \widehat{h}^+ + \widehat{h}^-$. 
\item[(2)]
$\widehat{h}^+\circ f = \delta \;\widehat{h}^+$, 
and $\widehat{h}^-\circ f^{-1} = \delta_- \;\widehat{h}^-$. 
\item[(3)]
$\widehat{h}^+ \geq 0$ and $\widehat{h}^- \geq 0$. 
\item[(4)] 
For $x \in \Aff^n(\overline{K})$, 
$\widehat{h}^+(x) =0$ if and only if 
$\widehat{h}^-(x) =0$ if and only if 
$\widehat{h}(x) =0$ if and only if 
$x$ is $f$-periodic.
\end{enumerate}
\end{Lemma}

\Proof
By the property (ii) in Theorem~\ref{thm:general}, 
we readily see (1). Let us see (2). 
By the property (ii), we have 
$\delta_- \widehat{h}(f^2(x)) + \delta \widehat{h}(x) = 
(1 + \delta \delta_-) \widehat{h}(f(x))$
and 
$\left(\frac{1}{\delta_-} + \delta\right) \widehat{h}(x)
= \widehat{h}(f(x)) + \frac{\delta}{\delta_-} \widehat{h}(f^{-1}(x))$ 
Taking the difference, we have 
\[
\delta_- \widehat{h}(f^2(x)) - \frac{1}{\delta_-}\widehat{h}(x)
= 
\delta \left(
\delta_- \widehat{h}(f(x)) - \frac{1}{\delta_-}\widehat{h}(f^{-1}(x))
\right). 
\]
This shows $\widehat{h}^+(f(x)) = \delta \;\widehat{h}^+(x)$. 
Similarly we have $\widehat{h}^+(f^{-1}(x)) 
= \delta_- \;\widehat{h}^-(x)$. 
Next let us see (3). 
Since $\widehat{h} \geq 0$ by Theorem~\ref{thm:general}, 
we have 
$\widehat{h}^+(f^l(x)) 
+ \widehat{h}^-(f^l(x))
= \widehat{h}(f^l(x)) \geq 0$ 
for any $l \in \ZZ$ and $x \in \Aff^n(\overline{K})$. 
This is equivalent to 
\[
\widehat{h}^+(x) \geq - \frac{1}{(\delta \delta_-)^l} 
\widehat{h}^-(x). 
\]
By letting $l \to  +\infty$, we have 
$\widehat{h}^+(x) \geq 0$. 
Similarly we have $\widehat{h}^-(x) \geq 0$. 

Next we will show (4). 
The assertion that ``$\widehat{h}(x) =0$ if and only if 
$x$ is $f$-periodic'' is shown in  Theorem~\ref{thm:general}. 
Since $\widehat{h}^+ \geq 0$ and $\widehat{h}^- \geq 0$, 
$0 = \widehat{h}(x) = \widehat{h}^+(x) + \widehat{h}^-(x)$ 
implies $\widehat{h}^+(x) =0$ and $\widehat{h}^-(x) =0$. 
We will see that $\widehat{h}^+(x)=0$ implies $\widehat{h}(x)=0$. 
A key observation here is that $\widehat{h}$ satisfies 
Northcott's finiteness property, which is a consequence of 
the property (i) of $\widehat{h}$ in Theorem~\ref{thm:general}. 
Suppose $\widehat{h}^+(x) =0$. Then 
\[
\widehat{h}(f^l(x)) = \widehat{h}^+(f^l(x)) 
+ \widehat{h}^-(f^l(x))
= \delta^l \widehat{h}^+(x) + 
\frac{1}{\delta_-^l} \widehat{h}^-(x)
= \frac{1}{\delta_-^l} \widehat{h}^-(x). 
\]
Let $L$ be a finite extension of $K$ over which 
$x$ is defined. Then 
\[
\{f^l(x) \in \Aff^n(\overline{K}) \mid l \geq 0\}
\quad\subseteq\quad  
\{y \in \Aff^n(\overline{K}) \mid 
\widehat{h}(y) \leq \widehat{h}^-(x)\}
\]
is finite. Hence $x$ is $f$-periodic. 
Similarly we see that $\widehat{h}^-(x)= 0$ implies 
$\widehat{h}(x)=0$.
\QED

For $x \in \Aff^n(\overline{K})$, we define 
the {\em $f$-orbit} of $x$ to be  
\[
O_f(x) := \{f^l(x) \mid l \in \ZZ\}. 
\]
Note that $O_f(x)$ is a finite set  
if and only if $x$ is $f$-periodic. 

For an $f$-orbit $O_f(x)$, we define the 
{\em canonical height of $O_f(x)$} to be 
\[
\widehat{h}(O_f(x)) 
= \frac{\log \widehat{h}^+(y)}{\log \delta}
+ \frac{\log \widehat{h}^-(y)}{\log \delta_-}
\qquad \in \RR \cup \{-\infty\}
\]
for any $y \in O_f(x)$. 

\begin{Lemma}
\label{lemma:height:of:orbit}
\begin{enumerate}
\item[(1)]
$\widehat{h}(O_f(x))$ is well-defined, i.e., 
$\widehat{h}(O_f(x))$ is independent of the choice 
of $y \in O_f(x)$. 
Moreover, $\widehat{h}(O_f(x)) = -\infty$ 
if and only if $O_f(x)$ is a finite set.
\item[(2)]
Assume $\# O_f(x) =  +\infty$. Then we have
\[
\widehat{h}(O_f(x)) + \epsilon_1 \leq 
\left(\frac{1}{\log \delta} + \frac{1}{\log \delta_-}\right)
\min_{y \in O_f(x)} \log \widehat{h}(y) \leq 
\widehat{h}(O_f(x)) + \epsilon_2, 
\]
where the constants $\epsilon_1$ and $\epsilon_2$ 
are given by 
\begin{align*}
\epsilon_1 & =
\frac{1}{\log \delta} \log\left(1 + \frac{\log \delta}{\log \delta_-}\right) 
+ \frac{1}{\log \delta_-} \log\left(1 + \frac{\log \delta_-}{\log \delta}\right), \\ 
\epsilon_2 & = \epsilon_1  
+ \left(\frac{1}{\log \delta} + \frac{1}{\log \delta_-} \right) 
\log\max\{\delta, \delta_-\}. 
\end{align*}
\end{enumerate}
\end{Lemma}

\Proof
(1) follows from Lemma~\ref{lemma:h:+:-}. 
To prove (2), set 
\[
p = 1 + \frac{\log \delta}{\log \delta_-}
\quad\text{and}\quad
q = 1 + \frac{\log \delta_-}{\log \delta}. 
\]
Then $p > 1$, $q >1$, and $\frac{1}{p} + \frac{1}{q} = 1$. 
Then we have 
\[
\widehat{h}(y) 
= \widehat{h}^+(y) + \widehat{h}^-(y)
= \frac{1}{p} \left(p^{\frac{1}{p}} \widehat{h}^+(y)^{\frac{1}{p}}\right)^p
+ \frac{1}{q} \left(q^{\frac{1}{q}} \widehat{h}^-(y)^{\frac{1}{q}}\right)^q 
\geq p^{\frac{1}{p}} q^{\frac{1}{q}} 
\widehat{h}^+(y)^{\frac{1}{p}}\widehat{h}^-(y)^{\frac{1}{q}}. 
\]
Hence, $\frac{1}{p} \log p + \frac{1}{q} \log q +  
\frac{1}{p}\log \widehat{h}^+(y) + \frac{1}{q}\log \widehat{h}^-(y)
\leq \log \widehat{h}(y)$. 
Since 
\[
\frac{1}{p} \log \widehat{h}^+(y) + \frac{1}{q} 
\log \widehat{h}^-(y) = 
\left(\frac{1}{\log \delta} 
+ \frac{1}{\log \delta_-}\right)^{-1} \widehat{h}(O_f(x)), 
\]
we obtain 
$\widehat{h}(O_f(x)) + \epsilon_1 \leq 
\left(\frac{1}{\log \delta} + \frac{1}{\log \delta_-}\right)
\min_{y \in O_f(x)} \log \widehat{h}(y)$.

On the other hand, we have 
$\widehat{h}(f^l(x)) 
= \delta^l \widehat{h}^+(x) + \delta_-^{-l} \widehat{h}^-(x)$ 
for $l \in \ZZ$. We set 
$g(t) = \delta^t \widehat{h}^+(x) + \delta_-^{-t} \widehat{h}^-(x)$
for $t \in \RR$, and 
\[
t_0 := 
\frac{\log(\widehat{h}^-(x)\log \delta_-) - 
\log(\widehat{h}^+(x) \log \delta)}{\log \delta + \log \delta_-}.
\]
Then one sees that $g$ takes its minimum at 
$t_0$, with $g(t_0) = p^{\frac{1}{p}} q^{\frac{1}{q}} 
\widehat{h}^+(x)^{\frac{1}{p}}\widehat{h}^-(x)^{\frac{1}{q}}$. 
Consequently as a function of $l \in \ZZ$, 
$\widehat{h}(f^l(x))$ takes its minimum 
at $l = [t_0]$ or $l = [t_0] +1$, where $[t_0]$ 
denotes the largest integer less than or equal to $t_0$. 
Then we get 
\begin{align*}
\widehat{h}(f^{[t_0]}(x))
& = \delta^{[t_0]} \widehat{h}^+(x) + \delta_-^{-[t_0]} \widehat{h}^-(x)
= \delta^{-(t_0 - [t_0])} \delta^{t_0} \widehat{h}^+(x) + 
\delta_-^{t_0 -[t_0]} \delta_-^{-t_0} \widehat{h}^-(x) \\
& < 
\max\{\delta, \delta_-\} 
\left(\delta^{t_0} \widehat{h}^+(x) + \delta_-^{-t_0} \widehat{h}^-(x) \right)
= \max\{\delta, \delta_-\} p^{\frac{1}{p}} q^{\frac{1}{q}} 
\widehat{h}^+(x)^{\frac{1}{p}}\widehat{h}^-(x)^{\frac{1}{q}}. 
\end{align*}
Similarly we get 
\begin{align*}
\widehat{h}(f^{[t_0]+1}(x))
& = \delta^{1 + [t_0] - t_0} \delta^{t_0} \widehat{h}^+(x) + 
\delta^{- (1 + [t_0] - t_0)} \delta_-^{-t_0} \widehat{h}^-(x) \\
& < \max\{\delta, \delta_-\} 
p^{\frac{1}{p}} q^{\frac{1}{q}} 
\widehat{h}^+(x)^{\frac{1}{p}}\widehat{h}^-(x)^{\frac{1}{q}}. 
\end{align*}
This shows $\left(\frac{1}{\log \delta} + \frac{1}{\log \delta_-}\right)
\min_{y \in O_f(x)} \log \widehat{h}(y)
\leq \widehat{h}(O_f(x)) + \epsilon_2$.
\QED

\begin{Theorem}
\label{thm:estimate:general}
Let $f: \Aff^n \to \Aff^n$ be a polynomial automorphism 
over a number field $K$ satisfying 
the conditions in Theorem~\ref{thm:general}, 
and $\widehat{h}$ a height 
function constructed in Theorem~\ref{thm:general}. 
Let $x$ be an element of $\Aff^n(\overline{K})$ 
such that $\# O_f(x) =  +\infty$. 
Then we have the following.
\begin{enumerate}
\item[(1)]
If $\left(\frac{1}{\log \delta} + \frac{1}{\log \delta_-} \right)\log T 
\geq \widehat{h}(O_f(x))$, then 
\[
\left\vert
\#\{ y \in O_f(x) \mid \widehat{h}(y) \leq T\} 
- \left(\frac{1}{\log \delta} + \frac{1}{\log \delta_-} \right)\log T
+ \widehat{h}(O_f(x))
\right\vert
\leq \frac{\log 2}{\log \delta} + \frac{\log 2}{\log \delta_-} + 1.
\]
Note that if $\left(\frac{1}{\log \delta} + \frac{1}{\log \delta_-} \right)\log T 
\leq \widehat{h}(O_f(x))$, it follows from 
Lemma~\ref{lemma:height:of:orbit}(2) that 
$\#\{ y \in O_f(x) \mid \widehat{h}(y) \leq T\} = \emptyset$. 
\item[(2)] 
$\displaystyle{
\#\{ y \in O_f(x) \mid h_{nv}(y) \leq T\}
= \left(\frac{1}{\log \delta} + \frac{1}{\log \delta_-} \right)\log T 
- \widehat{h}(O_f(x)) + O(1)}$
as $T\to +\infty$, 
where the $O(1)$ constant depends only on $f$ and the choice 
of $\widehat{h}$. 
\end{enumerate}
\end{Theorem}

\Proof
Since $\# O_f(x) =  +\infty$, the map 
$\ZZ \ni l \mapsto f^l(x) \in \Aff^n(\overline{K})$ 
is one-to-one. 
Then 
\begin{align*}
\#\{ y \in O_f(x) \mid \widehat{h}(y) \leq T\}
& = \#\{ l \in \ZZ \mid \widehat{h}(f^l(x)) \leq T\} \\
& = \#\{ l \in \ZZ \mid 
\delta^l \widehat{h}^+(x) + \delta_-^{-l} \widehat{h}^-(x) \leq T\}. 
\end{align*}
Then it follows from 
Lemma~\ref{lemma:number} that 
\[
-1 + \frac{\log{\frac{T}{2 \widehat{h}^+(x)}}}{\log \delta}
+  \frac{\log{\frac{T}{2 \widehat{h}^-(x)}}}{\log \delta_-}
\leq 
\#\{ y \in O_f(x) \mid \widehat{h}(y) \leq T\}
\leq 
1 + \frac{\log{\frac{T}{\widehat{h}^+(x)}}}{\log \delta}
+ \frac{\log{\frac{T}{\widehat{h}^-(x)}}}{\log \delta_-},   
\]
for $T \geq {\widehat{h}^+(x)}^{\frac{\log \delta_-}{\log \delta + \log \delta_-}} 
{\widehat{h}^-(x)}^{\frac{\log \delta}{\log \delta + \log \delta_-}}$ 
or equivalently $\left(\frac{1}{\log \delta} + \frac{1}{\log \delta_-} 
\right)\log T \geq \widehat{h}(O_f(x))$. 

On the other hand, we have 
\begin{gather*}
-1 + \frac{\log{\frac{T}{2 \widehat{h}^+(x)}}}{\log \delta}
+  \frac{\log{\frac{T}{2 \widehat{h}^-(x)}}}{\log \delta_-}
= -1 - \frac{\log 2}{\log \delta} - \frac{\log 2}{\log \delta_-} 
+ \left(\frac{1}{\log \delta} + \frac{1}{\log \delta_-} \right)\log T
- \widehat{h}(O_f(x)), \\ 
1 + \frac{\log{\frac{T}{\widehat{h}^+(x)}}}{\log \delta}
+ \frac{\log{\frac{T}{\widehat{h}^-(x)}}}{\log \delta_-}
= 1 + \left(\frac{1}{\log \delta} + \frac{1}{\log \delta_-} \right)\log T
- \widehat{h}(O_f(x)). 
\end{gather*}
Thus we obtain (1). Next, we will show (2). 
Since $h_{nv} \gg\ll \widehat{h}$ by the property (i) of 
Theorem~\ref{thm:main}, there exist a positive constant $a_2$ and 
a constant $b_2$ such that 
$\widehat{h} \leq a_2 h_{nv} + b_2$. 
Then we have 
\begin{align*}
& \#\{ y \in O_f(x) \mid h_{nv}(y) \leq T\} \\
&\quad \leq \#\{ y \in O_f(x) \mid  \widehat{h}(y) \leq a_2 T +b_2\} \\
&\quad \leq \left(\frac{1}{\log \delta} + \frac{1}{\log \delta_-} 
\right)\log (a_2T+b_2) 
- \widehat{h}(O_f(x)) + 1 + \frac{\log 2}{\log \delta} + \frac{\log 2}{\log \delta_-}\\
&\quad \leq 
\left(\frac{1}{\log \delta} + \frac{1}{\log \delta_-} \right)\log T
- \widehat{h}(O_f(x)) + O(1) 
\qquad \text{as $T\to +\infty$.}
\end{align*}
Using $a_1 h_{nv} + b_1 \leq \widehat{h}$ for some 
positive constant $a_1$ and constant $b_1$, we have 
$\#\{ y \in O_f(x) \mid h_{nv}(y) \leq T\}
\geq \left(\frac{1}{\log \delta} + \frac{1}{\log \delta_-} \right)\log T
- \widehat{h}(O_f(x)) + O(1)$ as $T\to +\infty$.
\QED

\begin{Lemma}
\label{lemma:number}
Let $A, B, T >0$ be positive numbers. 
If $T \geq A^{\frac{\log \delta_-}{\log \delta + \log \delta_-}} 
B^{\frac{\log \delta}{\log \delta + \log \delta_-}}$, 
then we have 
\[
-1 + \frac{\log{\frac{T}{2A}}}{\log \delta}
+  \frac{\log{\frac{T}{2B}}}{\log \delta_-}
\leq 
\#\{ l \in \ZZ \mid 
\delta^l A + \delta_-^{-l} B \leq T\}
\leq 
1 + \frac{\log{\frac{T}{A}}}{\log \delta}
+ \frac{\log{\frac{T}{B}}}{\log \delta_-}. 
\] 
\end{Lemma}

\Proof
If $l \in \ZZ$ satisfies $\delta^l A + \delta_-^{-l} B \leq T$, 
then $\delta^l A \leq T$ and $\delta_-^{-l} B \leq T$. 
Note that $\frac{\log{\frac{B}{T}}}{\log \delta_-} 
\leq \frac{\log{\frac{T}{A}}}{\log \delta}$ is equivalent to 
$T \geq A^{\frac{\log \delta_-}{\log \delta + \log \delta_-}} 
B^{\frac{\log \delta}{\log \delta + \log \delta_-}}$. 
Then, for $T \geq A^{\frac{\log \delta_-}{\log \delta + \log \delta_-}} 
B^{\frac{\log \delta}{\log \delta + \log \delta_-}}$, we have 
\[
\#\{ l \in \ZZ \mid \delta^l A + \delta_-^{-l} B \leq T\}
\leq 
\#\left\{ l \in \ZZ \;\left\vert\; \frac{\log{\frac{B}{T}}}{\log \delta_-} 
\leq l \leq \frac{\log{\frac{T}{A}}}{\log \delta} \right.\right\}
\leq 
1 + \frac{\log{\frac{T}{A}}}{\log \delta}
+ \frac{\log{\frac{T}{B}}}{\log \delta_-}.  
\]

On the other hand, if $l\in \ZZ$ satisfies 
$\delta^l A \leq \frac{T}{2}$ and $\delta_-^{-l} B \leq \frac{T}{2}$, 
then $\delta^l A + \delta_-^{-l} B \leq T$. Thus, 
\[
\#\{ l \in \ZZ \mid \delta^l A + \delta_-^{-l} B \leq T\}
\geq 
\#\left\{ l \in \ZZ \;\left\vert\; \frac{\log{\frac{2B}{T}}}{\log \delta_-} 
\leq l \leq \frac{\log{\frac{T}{2A}}}{\log \delta} \right. \right\}
\geq 
-1 + \frac{\log{\frac{T}{2A}}}{\log \delta}
+ \frac{\log{\frac{T}{2B}}}{\log \delta_-}. 
\]
\QED

\medskip
{\sl Proof of Theorem~\ref{thm:estimate}. }\quad
As we saw in the proof of 
Theorem~\ref{thm:main} and Corollary~\ref{cor:main}, 
polynomial automorphisms on $\Aff^2$ of dynamical 
degree $\geq 2$ satisfy the conditions 
in Theorem~\ref{thm:general}. 
Then Theorem~\ref{thm:estimate} follows from 
Theorem~\ref{thm:estimate:general}. 
\QED

\end{document}